\documentclass[12pt]{article}
\usepackage{amsmath}
\usepackage{amssymb}
\usepackage{amsthm}
\usepackage{amsfonts}
\usepackage{amscd}
\usepackage[mathscr]{eucal}
\newcommand{\Z} {{\mathbb  Z}}
\newcommand{\Q}{{\mathbb  Q}}
\newcommand{\F}{{\mathbb  F}}
\newcommand{\C}{{\mathbb  C}}
\newcommand{\N}{{\mathbb  N}}
\newcommand{\R} {{\mathbb R}}

\begin{document}
\parindent  25pt
\baselineskip  10mm
\textwidth  15cm    \textheight  23cm \evensidemargin -0.06cm
\oddsidemargin -0.01cm

\title{ { Geometry of Non-Archimedean Gromov-Hausdorff distance  }}
\author{\mbox{}
{ Derong Qiu }
\thanks{ \quad E-mail:
derong@mail.cnu.edu.cn } \\
(Department of Mathematics, Capital Normal University, \\
 Beijing 100037, P.R.China)  }

\date{}
\maketitle
\parindent  24pt
\baselineskip  10mm
\parskip  0pt

\par     \vskip  0.6 cm

\hspace{-0.6cm}{\bf 1. \ Introduction and statement of main
results }

\par \vskip 0.8 cm

In this paper, we study the geometry of non-Archimedean
Gromov-Hausdorff metric. This is the first part of our series
work, which we try to establish some facts about the counterpart
of Gromov-Hausdorff metric in the non-Archimedean spaces. One of
the motivation of this work is to find some implied relations
between this geometry and number theory via p-adic analysis, so
that we can use the former as a tool to study the relating
arithmetic aspects.
\par  \vskip 0.2 cm

Now we state the main results in the present work. Firstly, in
order to establish a compatible theory for the counterpart of
Gromov-Hausdorff metric in the non-Archimedean spaces, like the
well known case for general metric spaces ( see [G] and [BBI]), we
construct three corresponding key tools $ strong \ correspondence
$ ( see Def.2.11 ), \ $ strong \ \varepsilon -$isometry ( see
Def.2.22 ) and $ strong \ \varepsilon -approximations $ ( see
Def.3.4 ) for the non-Archimedean Gromov-Hausdorff distance $
\widehat{d}_{GH} $ ( see Def.2.1). Then we obtain the following
Theorems and formulae ( see also the related Cor.2.24 and 3.7 )
for explicitly computing $ \widehat{d}_{GH}. $
\par  \vskip 0.2 cm

{\bf Theorem A } ( see Theorem 2.14 below ). \ For any two
non-Archimedean metric spaces $ ( X , \widehat{d}_{X} ) $ and $ (
Y , \widehat{d}_{Y} ), $ \ $$ \widehat{d}_{GH} ( X , Y ) = \inf_{
\widehat{\textit{C }} }
 ( \text{dis} \widehat{ \textit{C} } ), $$  where the $ \inf $ is
 taken over all $ strong \ correspondences \ \widehat{ \textit{C }} $
 between $ X $ and $ Y. $
\par \vskip 0.2 cm

{\bf Theorem B } ( see Theorem 2.23 below ). \ Let $ ( X ,
\widehat{d}_{X} ) $ and $ ( Y , \widehat{d}_{Y} ) $ be two
non-Archimedean metric spaces and $ \varepsilon > 0. $ Then
\par  \vskip 0.15 cm
(1) \ If $ \widehat{d}_{GH} ( X , Y ) < \varepsilon , $ then there
exists a $ strong \ \varepsilon -$isometry from $ X $ to $ Y. $
\par  \vskip 0.15 cm
(2) \ If there exists a $ strong \ \varepsilon -$isometry from $ X
$ to $ Y, $ then $ \widehat{d}_{GH} ( X , Y ) \leq \varepsilon . $
\par  \vskip 0.2 cm

{\bf Theorem C } ( see Theorem 3.5 below ). \ Let $ ( X ,
\widehat{d}_{X} ) $ and $ ( Y , \widehat{d}_{Y} ) $ be two compact
non-Archimedean metric spaces.
\par  \vskip 0.15 cm
(1) \ If $ X $ and $ Y $ are $ strong \ \varepsilon
-approximations $ of each other, then $ \widehat{d}_{GH} ( X , Y )
\leq \varepsilon . $
\par  \vskip 0.15 cm
(2) \ If $ \widehat{d}_{GH} ( X , Y ) < \varepsilon , $ then $ X $
and $ Y $ are $ strong \ \varepsilon -approximations $ of each
other.
\par  \vskip 0.2 cm

Next we establish several convergence theorems for the
non-Archimedean metric spaces under $ \widehat{d}_{GH} $ ( see
Theorems 3.8 , 3.9 and 3.12 below ), one of them is the following
Compactness Theorem about the $ strongly \ uniformly \ totally \
bounded $ class ( see Def.3.11 for its definition ) of compact
non-Archimedean metric spaces.
\par  \vskip 0.2 cm

{\bf Theorem D } ( Compactness Theorem ) ( see Theorem 3.12 below ). \\
Any $ strongly \ uniformly \ totally \ bounded $ class $
\widehat{\mathcal{X}}_{sut} $ of compact non-Archimedean metric
spaces is pre-compact in the strong Gromov-Hausdorff topology.
That is, any sequence of elements of $ \widehat{\mathcal{X}}_{sut}
$ contains a Cauchy subsequence under the metric $
\widehat{d}_{GH}. $
\par  \vskip 0.2 cm

Using the above three tools, especially the $ strong \ \varepsilon
-$isometry and its related results, we obtain the following
theorem, which enables us conveniently calculate the
non-Archimedean Gromov-Hausdorff distance between non-Archimedean
metric spaces. As an application, we use it working out such
results about local fields ( see Example 4.7 below )
\par  \vskip 0.2 cm

{\bf Theorem E } ( see Theorem 4.2 below ). \ Let $ ( X, \
\widehat{d}_{X}) $ and $ ( Y, \ \widehat{d}_{Y}) $ be two
non-Archimedean metric spaces and denote $ D = \max \{
\text{diam}(X), \ \text{diam} (Y) \}. $ Then
\par  \vskip 0.15 cm
(1) \ $ \widehat{d}_{GH} ( X, \ Y ) \geq \inf \{ \varepsilon > 0 :
\ W_{X} ( X )_{ \geq \varepsilon } = W_{Y} ( Y )_{ \geq
\varepsilon } \} . $
\par  \vskip 0.15 cm
(2) \ If there does not exist $ \varepsilon > 0 $ such that $ \
W_{X} ( X )_{ \geq \varepsilon } =
 W_{Y} ( Y )_{ \geq \varepsilon } , $
then $ \widehat{d}_{GH} ( X, \ Y ) = \infty . $
\par  \vskip 0.15 cm
(3A) \ If $ D < + \infty , $ then $  \widehat{d}_{GH} ( X, \ Y )
\leq D. $
\par  \vskip 0.15 cm
(3B) \ If $ D < + \infty $ and $ \text{diam} (X) \neq \text{diam}
(Y) , $ then $ \widehat{d}_{GH} ( X, \ Y ) = D . $
\par  \vskip 0.2 cm
By Theorem E, we obtain the following theorem about converging
sequence.
\par  \vskip 0.2 cm

{\bf Theorem F }( see Theorem 4.4 and Remark 4.5 below ). \ Let $
\{ X_{n} \}_{n = 1 }^{\infty } $ be a sequence of non-Archimedean
metric spaces with $ \text{diam} ( X_{n} ) < + \infty $ for each $
n \in \N , $ and $ X $ be a non-Archimedean metric space with $
\text{diam} ( X ) < + \infty . $ If $ X_{n} \longrightarrow
_{\text{GH}_{S}} X, $ then
\par  \vskip 0.15 cm
(1) \ If $ \text{diam} ( X ) = 0 , $ then $ \text{diam} ( X_{n} )
\longrightarrow 0 $ as $ n \longrightarrow \infty . $
\par  \vskip 0.15 cm
(2) \ If $ \text{diam} ( X ) > 0 , $ then there exists a $ n_{0}
\in \N $ such that for all $ n > n_{0}, \ \text{diam} ( X_{n} ) =
\text{diam} ( X) . $
\par  \vskip 0.2 cm

Since we already have two kinds of metric structures $
\widehat{d}_{GH} $ and $d_{GH} $ for non-Archimedean metric
spaces, to compare them, we define their ratio as the metric ratio
function $ \Upsilon _{m} $ (see Def.2.16 below) and then by using
the explicit formula of Theorem E we prove that $ \Upsilon _{m} $
is unbounded, which is stated in the following theorem. This
Theorem is key to our construction, e.g., from it we know that the
strong convergence $ X_{n} \longrightarrow _{ \text{GH}_{S}} X $
and the convergence $ X_{n} \longrightarrow _{ \text{GH}} X $ (see
Def.3.1 below for $ \longrightarrow _{ \text{GH}_{S}} $  and [BBI]
for $ \longrightarrow _{ \text{GH}} $ ) are not equivalent, so the
two kinds of metric structures $ \widehat{d}_{GH} $ and $d_{GH} $
are different in essential.

{\bf Theorem G } ( see Theorem 4.8 below ). \ The metric ratio
function $ \Upsilon _{m} $ is unbounded, in other words, for any $
c \geq 2 , $ there exist non-Archimedean metric spaces $ X $ and $
Y $ such that $ \widehat{d}_{GH } ( X , Y ) \geq c \cdot d_{GH } (
X , Y ) . $
\par  \vskip 0.2 cm

In fact, in the proof of Theorem G, we obtain that $$ \Upsilon
_{m} ( \Z_{p}, \ \Z_{q}^{\Delta } ) = 2 q + 2 \longrightarrow
\infty \ \text{as} \ q \longrightarrow \infty  $$ and $ d_{GH} (
\Z_{p}, \ \Z_{q}^{\Delta } ) = \frac{1}{2 q } $ ( see the proof of
Theorem 4.8 and Remark 4.9 below ).
\par  \vskip 0.2 cm

Moreover, some questions about the relations between $
\widehat{d}_{GH} $ and $ d_{GH} $ and computation of $
\widehat{d}_{GH} $ are suggested ( see Questions 2.18 and 4.3 ). A
equilibrium function (see Def.2.19 and Lemma 2.20) is defined
associating to each $ stong \ correspondence \
\widehat{\mathcal{C}} $ between non-Archimedean metric spaces, and
a question about this function is also suggested (see Question
2.21).
\par  \vskip 0.2 cm

Based on these results, we will discuss the related arithmetic
applications as well as geometric structures on $ p-$adic
manifolds in separate papers [Q1] and [Q2].
\par \vskip 0.2 cm

{\bf Notation and terminology.} As usual, the symbols $ \Z , \Q ,
\R , \C , \F_{p} , \Z_{p} $ and $ \Q_{p} $ represent the integers,
rational numbers, real numbers, complex numbers, field with $ p $
elements, $ p-$adic integers and $ p-$adic numbers, respectively.
We denote the completion of the algebraic closure $ \overline{
\Q_{p}} $ of $ \Q_{p} $ by $ \C_{p}, $ which is endowed the
normalized valuation $ |.|_{p} $ satisfying $ |p |_{p} =
\frac{1}{p}, $ and called the Tate field ( see [K], [Se]). \\
Let $ X $ be a metric space endowed with the metric $ d . $ For
any subsets $ A $ and $ B $ of $ X , $ the diameter of $ A $ is
diam$(A) = \sup \{ d ( x, y ) : x, y \in A \}, \ A $ is bounded if
diam$( A ) < \infty . $ The distance between $ A $ and $ B $ is
dist$( A, B ) = \inf \{ d ( x, y ) : x \in A , y \in B \} , $ in
particular, for $ x \in X , $ dist$ ( x, A ) = $ dist$ ( \{ x \} ,
A ) . $ \ For a set $ Y \subset X $ and $ \varepsilon > 0 , $ its
$ \varepsilon -$neighborhood is the set $ U_{\varepsilon} ( Y ) =
\{ x \in X : \text{dist} ( x , Y ) < \varepsilon \}. $ \ $ Y $ is
an $ \varepsilon -$net in $ X $ if dist$( x , Y ) < \varepsilon $
for all $ x \in X . $  Particularly, $ X $ is called totally
bounded if for any $ \varepsilon > 0 $ there is a finite $
\varepsilon -$net in it. The Hausdorff distance between $ A $ and
$ B , $ denoted by $ d_{ H } ( A , B ) , $ is defined by
\begin{align*}
 d_{ H } ( A , B ) &= \inf \{ \varepsilon > 0 : A \subset U_{\varepsilon} ( B ) \
\text{ and } \ B \subset U_{\varepsilon} ( A ) \} \\
&= \max \{ \sup_{ a \in A } \text{dist} ( a , B ) , \quad  \sup_{
b \in B } \text{dist}( b , A ) \} .
\end{align*}
Given two metric spaces $ ( X , d_{X} ) $ and $ ( Y , d_{Y} ) , $
the Gromov-Hausdorff distance between them is defined as $$ d_{GH}
( X , Y ) = \inf \{ d_{Z, H} ( f ( X ), g( Y ) ) \} , $$ where the
$ \inf $ is taken over all metric spaces $ ( Z , d_{Z} ) $ and
over all isometric embeddings $ f : X \hookrightarrow Z $ and $ g
: Y \hookrightarrow Z . \ d_{Z, H} $ denotes the Hausdorff
distance between subsets of $ Z $ ( see [BBI] and [G]).
\par \vskip 0.2 cm

Throughout this paper, unless otherwise specified, we use the
notation $ \widehat{d} $ denote the non-Archimedean metric on a
space $ X , $ that is, $ \widehat{d} $ is a metric on $ X $ and
satisfies the strong triangle inequality $$  \widehat{d} ( x , z )
\leq \max \{ \widehat{d} ( x , y ) , \  \widehat{d} ( y , z ) \} \
( x , y , z \in X ) . $$ A set $ X $ endowed with a
non-Archimedean metric $ \widehat{d}_{X} $ is called a
non-Archimedean metric space ( also called ultrametric space ),
which is denoted by $ ( X , \widehat{d}_{X} ) . $ \ For the basic
properties of non-Archimedean metric spaces, we refer to [BGR] and
[Sc].

\par  \vskip 0.2 cm

{\bf Acknowledgements.} \ This work began on the summer of the
year 2007, when Professor Xiaochun Rong (Rutgers Univ. USA)
visited our mathematics department. I thank him heartily for
helpful discussing on metric geometry, especially for letting me
know the Gromov geometry as well as the two books [G] and [BBI],
from which I benefit very much. At that time, I started to be
interested in the question on how to establish the counterpart of
Gromov-Hausdorff metric in the non-Archimedean spaces, so that we
can use it as a tool to study arithmetic aspects in number theory,
particularly in elliptic curve (see [Q1] for our work on the part
of arithmetic applications), which is one of my main interesting
research fields. I thank heartily Prof. Xianke Zhang (Tsinghua
Univ.), my PhD Advisor, for introducing me in the field of
algebraic number theory, and for his great help and encouragement.
I thank heartily Prof. Qingzhong Li, Prof. Kezheng Li and Prof. Ke
Wu (Capital Normal Univ., Beijing) for their warm-hearted help and
regard. I thank heartily Prof. Shouwu Zhang (Columbia Univ. of New
York) and Prof. Keqin Feng (Tsinghua Univ.) for their generous
help and support in my studying arithmetic geometry at Columbia
Univ. of New York (2004-2005) and ICTP, Italy (2002).

\par     \vskip  1 cm

\hspace{-0.6cm}{\bf 2. \ Non-Archimedean Gromov-Hausdorff
distance }

\par \vskip 0.8 cm

Let $ ( X_{i} , \widehat{d}_{i} ) \ ( i = 1, 2 ) $ be two
non-Archimedean metric spaces. As in the general case of metric
spaces, it is not difficult to show that there exist a
non-Archimedean metric space $ ( X , \widehat{d}_{X})$ such that
both $ X_{1} $ and $ X_{2} $ can be isometrically embedded into $
X . $
\par \vskip 0.2 cm

{\bf Definition 2.1. } \ Let $ ( X_{i} , \widehat{d}_{i} ) \ ( i =
1, 2 ) $ be two non-Archimedean metric spaces. The non-Archimedean
Gromov-Hausdorff distance between them, denoted by $
\widehat{d}_{GH} ( X_{1} , X_{2} ) , $ is defined by
$$ \widehat{d}_{GH} ( X_{1} , X_{2} ) =
\inf \{ \widehat{d}_{H} ( f_{1}( X_{1} ), f_{2}( X_{2} ) ) \} , $$
where the $ \inf $ is taken over all non-Archimedean metric spaces
$ ( X , \widehat{d}) $ and over all isometric embeddings $ f_{1} :
X_{1} \hookrightarrow X $ and $ f_{2} : X_{2} \hookrightarrow X .
\ \widehat{d}_{H} $ denotes the Hausdorff distance between subsets
of $ X. $ \\
By definition, obviously we always have $ d_{GH} ( X_{1} , X_{2} )
\leq \widehat{d}_{GH} ( X_{1} , X_{2} ) $ for all non-Archimedean
metric spaces $ X_{1} $ and $ X_{2} . $
\par \vskip 0.2 cm

{\bf Example 2.2.} \ Let $ ( X , \widehat{d} )  $ be a
non-Archimedean metric space, $ Y \subset X $ be an $ \varepsilon
-$net with $ \varepsilon > 0 . $ Then easily we have $
\widehat{d}_{GH} ( X , Y ) \leq  \varepsilon . $
\par \vskip 0.2 cm

{\bf Definition 2.3. } \ Let $ ( X_{i} , \widehat{d}_{i} ) \ ( i =
1, 2 ) $ be two non-Archimedean metric spaces. Denote $ X = X_{1}
\sqcup X_{2} $ be the disjoint union. A non-Archimedean metric $
\widehat{d} $ on $ X $ is called admissible if it extends the
metrics on $ X_{1} $ and $ X_{2} , $ i.e., $ \widehat{d} \mid
_{X_{i}} = \widehat{d}_{i} \ ( i = 1, 2 ) . $ \ More generally,
for a family of non-Archimedean metric spaces $ \{ (X_{i} ,
\widehat{d}_{i} ) \}_{i \in \Lambda } ,  $ let $ X = \sqcup_{i \in
\Lambda } X_{i} $ be the disjoint union of all $ X_{i} . $ A
non-Archimedean metric $ \widehat{d} $ on $ X $ is called
admissible if it extends the metrics on all $ X_{i} , $ i.e., $
\widehat{d} \mid _{X_{i}} = \widehat{d}_{i} $ for every $ i \in
\Lambda . $ \\
Then as in the usual case [G], we define
\par \vskip 0.2 cm

{\bf Definition 2.4. } \ Let $ ( X_{i} , \widehat{d}_{i} ) \ ( i =
1, 2 ) $ be two non-Archimedean metric spaces. Define $$
\overline{\widehat{d}}_{GH}( X_{1} , X_{2} ) = \inf \{
\widehat{d}_{H} ( X_{1},  X_{2} ) \} ,
$$ where the $ \inf $ is taken over all admissible non-Archimedean
metrics $ \widehat{d} $ on $ X_{1} \sqcup X_{2} . $
\par \vskip 0.2 cm

{\bf Proposition 2.5. } \ Let $ ( X_{i} , \widehat{d}_{i} ) \ ( i
= 1, 2 ) $ be two non-Archimedean metric spaces. Then $$
\widehat{d}_{GH}( X_{1} , X_{2} ) = \overline{\widehat{d}}_{GH}(
X_{1} , X_{2} ) . $$  \\
{\bf Proof. } \ It is obvious that $ \widehat{d}_{GH}( X_{1} ,
X_{2} ) \leq \overline{\widehat{d}}_{GH}( X_{1} , X_{2} ) . $ So
we only need to prove the converse. For any $ \varepsilon > 0 , $
there exist a non-Archimedean metric space $ ( X , \widehat{d} ) $
and isometric embeddings $ f_{i} : X_{i} \rightarrow X \ ( i = 1,
2) $ such that $ \widehat{d}_{H}( f_{1} ( X_{1} ) , f_{2} ( X_{2}
) ) < \widehat{d}_{GH}( X_{1} , X_{2} ) + \varepsilon . $ Denote $
X^{\prime } = X_{1} \sqcup X_{2} $ ( disjoint union ), we define a
function $ \widehat{d_{3}} $ on $ X^{\prime } \times X^{\prime } $
as follows: For any $ x_{1} , x_{2} \in X^{\prime } , $  \
$$ \widehat{d_{3}}( x_{1}, x_{2}) = \widehat{d_{3}}( x_{2}, x_{1})
= \left\{
\begin{array}{l} \widehat{d} ( f_{i} ( x_{1} ) , f_{i} ( x_{2}) )
\quad \text{if} \
x_{1} , x_{2} \in X_{i} \ ( i = 1 , 2 ) ; \\
\max \{ \widehat{d} ( f_{1} ( x_{1} ) , f_{2} ( x_{2}) ) , \
\varepsilon \} \quad \text{if} \ x_{1} \in X_{1} , x_{2} \in X_{2}
.
\end{array}
\right. $$ Then by a direct calculation, it is not difficult to
show that $ \widehat{d_{3}} $ is an admissible non-Archimedean
metric on $ X^{\prime } . $ Therefore the Hausdorff distance
between $ X_{1} $ and $ X_{2} $ in $ X^{\prime } $ is
\begin{align*}
&\widehat{d_{3}}_{H}( X_{1}, X_{2}) = \max \{ \sup_{ x_{1} \in
X_{1} } \text{dist} ( x_{1} , X_{2} ) , \quad  \sup_{ x_{2} \in
X_{2}} \text{dist}( x_{2} , X_{1} ) \} \\
&\leq \max \{ \sup_{ x_{1} \in X_{1} } \max \{ \text{dist} ( f_{1}
( x_{1}) , f_{2} (X_{2} ) ) , \ \varepsilon \}, \ \sup_{ x_{2} \in
X_{2} } \max \{ \text{dist} ( f_{2} ( x_{2}) , f_{1} ( X_{1} ) ) ,
\ \varepsilon \}
\} \\
&= \max \{  \sup_{ x_{1} \in X_{1} } \text{dist} ( f_{1} ( x_{1})
, f_{2} (X_{2} ) ), \quad  \sup_{ x_{2} \in X_{2} } \text{dist} (
f_{2} ( x_{2}) ,
f_{1} ( X_{1} ) ) , \quad \varepsilon  \} \\
&= \max \{ \widehat{d}_{H} ( f_{1} ( X_{1} ), \ f_{2} (X_{2})),
\quad \varepsilon \} \leq \max \{  \widehat{d}_{GH} (X_{1} ,
X_{2} ) + \varepsilon , \ \varepsilon \} \\
&= \widehat{d}_{GH} (X_{1} , X_{2} ) + \varepsilon .
\end{align*}
Hence by definition, $ \overline{\widehat{d}}_{GH}( X_{1} , X_{2}
) \leq \widehat{d}_{GH} (X_{1} , X_{2} ) + \varepsilon , $ so $
\overline{\widehat{d}}_{GH}( X_{1} , X_{2} ) \leq \widehat{d}_{GH}
(X_{1} , X_{2} ) $ because $ \varepsilon > 0 $ is arbitrary. The
proof is completed. \quad $ \Box $
\par \vskip 0.2 cm

{\bf Proposition 2.6 } [Zarichnyi]. \ The function $
\widehat{d}_{GH} $ is a non-Archimedean metric on the set of
isometry classes of non-Archimedean metric spaces.
\par \vskip 0.2 cm

{\bf Proof. } \ This can be verified directly via a tedious
calculation as I have done myself. After finishing this paper, I
read a paper [Z] of I. Zarichnyi, in which he has already proven
this proposition ( see Theorem 1.2 of [Z]) as well as defined the
function $ \widehat{d}_{GH}, $ which he called the
Gromov-Hausdorff ultrametric. So this proposition owes completely
to Zarichnyi, and we refer to his paper [Z] for the detailed
proof. \quad $ \Box $
\par \vskip 0.2 cm

For any non-Archimedean metric spaces $ ( X_{i} , \widehat{d}_{i}
) \ ( i = 1, 2 ), $ they have now been defined two kinds of
distances $ d_{GH} ( X_{1} , X_{2} ) $ and $ \widehat{d}_{GH} (
X_{1} , X_{2} ), $ the former is their Gromov-Hausdorff distance,
and the later is their non-Archimedean Gromov-Hausdorff distance.
We will study the relations between them (see Theorem 2.8,
Corollary 2.9 and Theorem 4.8 below). In the following, we denote
by $ \widehat{\Gamma }_{c} $ the set of isometry classes of all
compact non-Archimedean metric spaces, and call $ (
\widehat{\Gamma }_{c} , \ \widehat{d}_{GH} ) $ (i.e., endowed with
the metric $ \widehat{d}_{GH} $) the non-Archimedean
Gromov-Hausdorff space. It is easy to see that $ \widehat{d}_{GH}
$ is a finite metric on this space. In fact, $ \widehat{d}_{GH}
(X_{1} , X_{2} ) < \infty $ for any bounded non-Archimedean metric
spaces $ ( X_{i} , \widehat{d}_{i} ) ( i = 1, 2 ) . $ Zarichnyi
has proved that $ ( \widehat{\Gamma }_{c} , \widehat{d}_{GH} ) $
is complete but not separable ( see [ Z, Prop. 2.1 and 2.2]).
\par \vskip 0.2 cm

{\bf Example 2.7. } \ Let $ ( X , \widehat{d} ) $ be a
non-Archimedean metric space, $ P $ be a non-Archimedean metric
space consisting of one point. Then $ \widehat{d}_{GH} ( X , P ) =
$ diam $ X , $ in particular, $ \widehat{d}_{GH} ( X , P ) = 2
d_{GH} ( X , P ). $
\par \vskip 0.2 cm

{\bf Proof. } \ The first equality can be easily verified by
Def.2.4 and Prop.2.5. The second equality follows from the fact
that $ d_{GH} ( X , P ) = \frac{1}{2} $ diam $ X $ ( see [BBI,
P.255]). \quad $ \Box $
\par \vskip 0.2 cm

Next we come to establish some formulae for explicitly computing
the non-Archimedean Gromov-Hausdorff distance. The first one is
via the key tool $ strong $  $ correspondence $ which we will
define in the following (see Def.2.11 below). Before doing this,
we first prove an inequality (see Theorem 2.8) which giving a
lower bound for $ \widehat{d}_{GH} $ by using the tool $
correspondence $ in the usual sense. For two sets $ X $ and $ Y ,
$ recall that [BBI, p.256, 257] a $ correpondence $ between them
is a set $ \textit{C} \subset X \times Y $ satisfying the
following condition: for every $ x \in X $ there exists at least
one $ y \in Y $ such that $ ( x , y ) \in \textit{C} , $ and
similarly for every $ y \in Y $ there exists an $ x \in X $ such
that $ ( x , y ) \in \textit{C} . $ For example, if $ f : X
\longrightarrow Y $ is a surjective map, then $ \textit{C}_{f} =
\{ ( x, f(x) ) : x \in X \} $ is a $ correspondence $ between $ X
$ and $ Y , $ which is called the $ correspondence $ associated
with $ f $ (see [BBI, p.256]). Now let $ \textit{C} $ be a $
correpondence $ between metric spaces $ ( X , d_{X} ) $ and $ ( Y
, d_{Y} ), $ then the $ distorsion $ of $ \textit{ C } $ is
defined by $$ \text{dis}\textit{C } = \sup \{ \mid d_{X} ( x ,
x^{\prime } ) - d_{Y} ( y , y^{\prime } ) \mid \ : \ ( x, y ),  (
x^{\prime } , y^{\prime } ) \in \textit{ C } \}. $$ Moreover, if $
f : X \longrightarrow Y $ is an arbitrary map, the $ distorsion $
of $ f $ is defined by $$ \text{dis}f = \sup_{x_{1} , x_{2} \in X
} \mid d_{Y} ( f(x_{1}) , f(x_{2}) - d_{X} ( x_{1} , x_{2} ) \mid
.
$$ It is easy to see that $ \text{dis}\textit{C}_{f} = \text{dis}f $
if $ f $ is a surjective map (See [BBI, p.257 and p.249] for the
properties of $ distorsions $ ).
\par \vskip 0.2 cm

{\bf Theorem 2.8. } \ For any two non-Archimedean metric spaces $
( X , \widehat{d}_{X} ) $ and $ ( Y , \widehat{d}_{Y} ), $ \  $$
 \widehat{d}_{GH} ( X , Y ) \geq \inf_{\textit{C} }
 ( \text{dis}\textit{C} ) , $$ where the $ \inf $ is taken over
 all $ correspondences $ \ $ \textit{C} $ between $ X $ and $ Y . $
\par \vskip 0.2 cm

{\bf Proof. } \ For any $ r > \widehat{d}_{GH} ( X , Y ) , $ by
definition, there exists a non-Archimedean metric spaces $ ( Z ,
\widehat{d}_{Z} ) $ such that $ X $ and $ Y $ can be isometrically
embedded in it and $ \widehat{d}_{Z, H} ( X , Y ) < r . $ We may
view $ X $ and $ Y $ as subspaces of $ Z $ with $ \widehat{d}_{Z}
\mid _{X} = \widehat{d}_{X} $ and $ \widehat{d}_{Z} \mid _{Y} =
\widehat{d}_{Y} . $ Let $ \textit{C}_{0} = \{ ( x , y ) :  x \in X
, \ y \in Y , \ \widehat{d}_{Z} ( x , y ) < r  \}. $ By $
\widehat{d}_{Z, H} ( X , Y ) < r $ we have $ X \subset U_{r} ( Y )
$ and $ Y \subset U_{r} ( X ) . $ So for any $ x \in X , $ there
exists a $ y \in Y $ such that $ \widehat{d}_{Z} ( x , y ) < r , $
hence $ ( x , y ) \in \textit{C}_{0} . $ Similarly, for any $ y
\in Y , $ there exists an $ x \in X $ such that $ ( x , y ) \in
\textit{C}_{0} . $ This shows that $ \textit{C}_{0} $ is a $
correspondence $ between $ X $ and $ Y . $ Now let $ ( x , y ) , (
x^{\prime }, y^{\prime } ) \in \textit{C}_{0} , $ then $
\widehat{d}_{Z} ( x , y ) < r $ and $ \widehat{d}_{Z} ( x^{\prime
}, y^{\prime } ) < r . $ So $ \widehat{d}_{X} ( x , x^{\prime } )
\leq \max \{ \widehat{d}_{Z} ( x , y ) , \ \widehat{d}_{Z} (
x^{\prime } , y ) \} \leq  \max \{ \widehat{d}_{Z} ( x , y ) , \
\widehat{d}_{Z} ( x^{\prime } , y^{\prime } ) , \  \widehat{d}_{Z}
( y^{\prime } , y ) \} \leq \max \{ \widehat{d}_{Z} ( x , y ) , \
\widehat{d}_{Z} ( x^{\prime } , y^{\prime }) \} + \widehat{d}_{Y}
( y^{\prime } , y ) , $ \ so $ \widehat{d}_{X}( x , x^{\prime }) -
\widehat{d}_{Y}( y , y^{\prime }) \\ \leq \max \{ \widehat{d}_{Z}(
x, y ) , \ \widehat{d}_{Z} ( x^{\prime }, y^{\prime })\} < r. $
Similarly, we have $ \widehat{d}_{Y} ( y , y^{\prime } ) -
\widehat{d}_{X} ( x , x^{\prime } ) < r. $ So $ \mid
\widehat{d}_{X} ( x , x^{\prime } ) - \widehat{d}_{Y}( y ,
y^{\prime }) \mid \ < r . $ Therefore by definition, $
\text{dis}\textit{C}_{0} = \sup \{ \mid \widehat{d}_{X} ( x ,
x^{\prime } ) - \widehat{d}_{Y} ( y , y^{\prime } ) \mid \ : \ ( x
, y) , \ ( x^{\prime } , y^{\prime } ) \in \textit{C}_{0} \} \leq
r . $ So $ r \geq \inf_{\textit{C} } ( \text{dis}\textit{C} ) . $
Since $ r > \widehat{d}_{GH} ( X , Y ) $ is arbitrary, we obtain $
\widehat{d}_{GH} ( X , Y ) \geq \inf_{\textit{C} } (
\text{dis}\textit{C} ). $ This proves Theorem 2.8.  \quad $ \Box $
\par \vskip 0.2 cm

It is well known that $  d_{GH} ( X , Y ) =
\frac{1}{2}\inf_{\textit{C} } ( \text{dis}\textit{C} )  $ ( see
[BBI, Thm.7.3.25]). So we have
\par \vskip 0.2 cm

{\bf Corollary 2.9. } \ For any two non-Archimedean metric spaces
$ ( X , \widehat{d}_{X} ) $ and $ ( Y , \widehat{d}_{Y} ), $ \  $$
 \widehat{d}_{GH} ( X , Y ) \geq 2 d_{GH} ( X , Y ). $$  Note that
Example 2.7 shows that the equality can hold in some cases.
\par \vskip 0.2 cm

Let $ X $ and $ Y $ be metric spaces and $ \varepsilon > 0 . $
Recall that a ( possibly noncontinuous ! ) map $ f : X
\longrightarrow  Y $ is called an $ \varepsilon -$isometry if
dis$f < \varepsilon $ and $ f ( X ) $ is an $ \varepsilon -$net in
$ Y $ ( see [BBI, p.258]. Note that we use $ < $ instead of $ \leq
$ in this definition).
\par \vskip 0.2 cm

{\bf Corollary 2.10. } \ Let $ ( X , \widehat{d}_{X} ) $ and $ ( Y
, \widehat{d}_{Y} ) $ be two non-Archimedean metric spaces and $
\varepsilon > 0 . $ If $ \widehat{d}_{GH} ( X , Y ) < \varepsilon
, $ then there exists an $ \varepsilon -$isometry from $ X $ to $
Y . $
\par \vskip 0.2 cm

{\bf Proof. } \ If $ \widehat{d}_{GH} ( X , Y ) < \varepsilon , $
then by Corollary 2.9, we have $ d_{GH} ( X , Y ) < \frac{1}{2}
\varepsilon  , $ so by [BBI, Cor.7.3.28], there exists an $
\varepsilon -$isometry from $ X $ to $ Y . $ \quad $ \Box $
\par \vskip 0.2 cm

{\bf Definition 2.11. } \ Let $ ( X , \widehat{d}_{X} ) $ and $ (
Y , \widehat{d}_{Y} ) $ be two non-Archimedean metric spaces, $
\textit{C} $ be a $ correspondence $ between $ X $ and $ Y . $ If
$ \textit{C} $ satisfies the following condition
\par \vskip 0.15 cm
( $ C_{NA} $ ) \ For any $ ( x , y ) \in X \times Y \setminus
 \textit{C} , $ there exist $ x^{\prime} \in X $ and  $ y^{\prime} \in Y
 $ such that $$ ( x , y^{\prime} ) , \  ( x^{\prime} , y )
 \in \textit{C} ,  \quad \text{and} \quad \widehat{d}_{X}
 ( x , x^{\prime} ) = \widehat{d}_{Y} ( y , y^{\prime} ) >
 \text{dis}\textit{C} . $$ Then we call $ \textit{C} $ a $ strong
\ correspondence $ between $ X $ and $ Y . $
\par \vskip 0.2 cm

{\bf Lemma 2.12. } \ The condition ( $ C_{NA} $ ) in Definition
2.11 is equivalent to the following condition
\par \vskip 0.15 cm
( $ C^{\prime}_{NA} $ ) \ If $ ( x , y ) \in X \times Y \setminus
\textit{C} , $ then for any $ x^{\prime} \in X $ and $ y^{\prime}
\in Y $ satisfying $ ( x , y^{\prime} ) , \  ( x^{\prime} , y )
\in \textit{C} , $ \ we have \ $  \widehat{d}_{X}
 ( x , x^{\prime} ) = \widehat{d}_{Y} ( y , y^{\prime} ) >
 \text{dis}\textit{C} . $
\par \vskip 0.2 cm

{\bf Proof. } \ ( $ C_{NA} ) \Longrightarrow  ( C^{\prime}_{NA} $
) : If $ ( x , y ) \in X \times Y \setminus \textit{C} , $ then
there exist $ x_{0} \in X $ and $ y_{0} \in Y $ such that $ ( x ,
y_{0} ) , \  ( x_{0} , y ) \in \textit{C} , $ \ and \ $
\widehat{d}_{X} ( x , x_{0} ) = \widehat{d}_{Y} ( y , y_{0} ) >
\text{dis}\textit{C} . $ Now for any $ x^{\prime} \in X $ and $
y^{\prime} \in Y $ satisfying $ ( x , y^{\prime} ) , \  (
x^{\prime} , y ) \in \textit{C} , $ \ we have \ $  \widehat{d}_{Y}
 ( y_{0} , y^{\prime} ) = \mid  \widehat{d}_{Y}
 ( y_{0} , y^{\prime} ) - \widehat{d}_{X}
 ( x , x ) \mid \leq \text{dis}\textit{C} . $ So by the strong
 triangle inequality, $  \widehat{d}_{Y} ( y , y^{\prime} ) =
 \max \{  \widehat{d}_{Y} ( y , y_{0} ) , \  \widehat{d}_{Y}
 ( y_{0} , y^{\prime} ) \} = \widehat{d}_{Y} ( y , y_{0} ) >
 \text{dis}\textit{C} . $ Similarly, $ \widehat{d}_{X} ( x , x^{\prime} )
 = \widehat{d}_{X} ( x , x_{0} ) > \text{dis}\textit{C} . $ Hence
$  \widehat{d}_{X} ( x , x^{\prime} ) = \widehat{d}_{Y} ( y ,
y^{\prime} ) > \text{dis}\textit{C} . $  \\
( $ C^{\prime}_{NA}) \Longrightarrow  ( C_{NA} $ ) : \ Obvious.
\quad $ \Box $
\par \vskip 0.2 cm

{\bf Notation. } \ In the following, we use the notation $
\widehat{ \textit{C} } $ instead of $ \textit{C} $ to denote a $
strong \ correspondence $ between $ X $ and $ Y . $
\par \vskip 0.2 cm

{\bf Lemma 2.13. } \ Let $ ( X , \widehat{d}_{X} ) $ and $ ( Y ,
\widehat{d}_{Y} ) $ be two non-Archimedean metric spaces.
\par \vskip 0.15 cm
(1) \ There always exist $ strong \ correspondences $ between $ X
$ and $ Y , $
\par \vskip 0.15 cm
(2) \ Let $ \widehat{ \textit{C} } $ be a $ strong \
correspondence $ between $ X $ and $ Y , $ let $ ( x, y ) , ( x,
y^{\prime} ) , ( x, y^{\prime\prime} )  \in X \times Y . $ If $ (
x, y^{\prime} ) , ( x, y^{\prime\prime} )  \in \widehat{
\textit{C} }  $ and $ ( x, y ) \notin \widehat{ \textit{C} } , $
then $ \widehat{d}_{Y} ( y , y^{\prime} ) = \widehat{d}_{Y} ( y ,
y^{\prime\prime} ) > \text{dis}\widehat{ \textit{C} } . $
Similarly, if $ ( x^{\prime}, y ) , ( x^{\prime\prime}, y )  \in
\widehat{ \textit{C} }  $ and $ ( x, y ) \notin \widehat{
\textit{C} } , $ then $ \widehat{d}_{X} ( x , x^{\prime} ) =
\widehat{d}_{X} ( x , x^{\prime\prime} ) > \text{dis}\widehat{
\textit{C} } . $
\par \vskip 0.2 cm

{\bf Proof. } \ ( 1 ) \ Obvious, e.g. $ \widehat{ \textit{C} } = X
\times Y . $
\par \vskip 0.15 cm
(2) \ Since $ ( x , y ) \in X \times Y \setminus \widehat{
\textit{C} }, $ by definition, there exist $ x_{0} \in X $ and $
y_{0} \in Y $ such that $ ( x, y_{0} ) , ( x_{0}, y ) \in
\widehat{ \textit{C} }  $ and $ \widehat{d}_{X} ( x , x_{0} ) =
\widehat{d}_{Y} ( y , y_{0} ) > \text{dis}\widehat{ \textit{C} } .
$ Then $  \widehat{d}_{Y} ( y_{0} , y^{\prime} ) = \ \mid
\widehat{d}_{Y} ( y_{0} , y^{\prime} ) - \widehat{d}_{X}
 ( x , x ) \mid \ \leq \text{dis}\widehat{ \textit{C} } . $
 So by the strong triangle inequality, we have
 $  \widehat{d}_{Y} ( y , y^{\prime} ) =
 \max \{  \widehat{d}_{Y} ( y , y_{0} ) , \  \widehat{d}_{Y}
 ( y_{0} , y^{\prime} ) \} = \widehat{d}_{Y} ( y , y_{0} ) >
 \text{dis}\widehat{ \textit{C} } . $ Also,
 $ \widehat{d}_{Y} ( y^{\prime} , y^{\prime\prime} )
 = \ \mid \widehat{d}_{Y} ( y^{\prime} , y^{\prime\prime} )
 - \widehat{d}_{X} ( x , x ) \mid \ \leq \text{dis}
 \widehat{ \textit{C} } . $
 So $ \widehat{d}_{Y} ( y , y^{\prime\prime} ) =
 \max \{ \widehat{d}_{Y} ( y , y^{\prime} ) , \  \widehat{d}_{Y}
 ( y^{\prime} , y^{\prime\prime} ) \} =
 \widehat{d}_{Y} ( y , y^{\prime} ) > \text{dis}\widehat{ \textit{C} } . $
\quad $ \Box $
\par \vskip 0.2 cm

{\bf Theorem 2.14. } \ For any two non-Archimedean metric spaces $
( X , \widehat{d}_{X} ) $ and $ ( Y , \widehat{d}_{Y} ), $ \ $$
\widehat{d}_{GH} ( X , Y ) = \inf_{ \widehat{\textit{C }} }
 ( \text{dis} \widehat{ \textit{C} } ), $$  where the $ \inf $ is taken over
 all $ strong \ correspondences \ \widehat{ \textit{C }} $
between $ X $ and $ Y . $
\par \vskip 0.2 cm

{\bf Proof. } \ Step 1. \ For any $ r > \widehat{d}_{GH} ( X , Y )
 , $ by definition, there exists a non-Archimedean metric space
 $ ( Z , \widehat{d}_{Z} ) $ such that $ X $ and $ Y $ can be
 isometrically embedded in $ Z $ ( and then we may view $ X $ and
 $ Y $ as subspaces of $ Z $ ) such that $ \widehat{d}_{Z, H} ( X , Y ) < r.
 $ Define $$ \widehat{\textit{C}}_{0} =
 \{ ( x , y ) \ : \ x \in X , \ y \in Y , \ \widehat{d}_{Z} ( X , Y ) \leq r  \}.
 $$ Then as done in the proof of Theorem 2.8, it is easy to verify
 that $ \widehat{\textit{C}}_{0} $ is a $ correspondence $ between
$ X $ and $ Y , $ and $ \text{dis} \widehat{ \textit{C}}_{0} \leq
r . $ Now let $ ( x , y ) \in X \times Y \setminus
\widehat{\textit{C}}_{0}. $ Then $ \widehat{d}_{Z} ( x , y ) >  r
. $ Since $ \widehat{\textit{C}}_{0} $ is a $ correspondence , $
there exist $ x^{\prime} \in X , $ and $ y^{\prime} \in Y $ such
that $ ( x^{\prime} , y ) , \ ( x , y^{\prime} ) \in \widehat{
\textit{C}}_{0} . $ So $ \widehat{d}_{Z} ( x^{\prime} , y ) \leq r
$ and $ \widehat{d}_{Z} ( x , y^{\prime} ) \leq  r . $ Then by
strong triangle inequality,
\begin{align*}
&\widehat{d}_{X} ( x, x^{\prime} ) = \widehat{d}_{Z} ( x ,
x^{\prime} ) = \max \{ \widehat{d}_{Z} ( x , y ), \
\widehat{d}_{Z} ( x^{\prime} , y ) \} = \widehat{d}_{Z} ( x , y )
> r .  \\
&\widehat{d}_{Y} ( y, y^{\prime} ) = \widehat{d}_{Z} ( y ,
y^{\prime} ) = \max \{ \widehat{d}_{Z} ( x , y ), \
\widehat{d}_{Z} ( x , y^{\prime} ) \} = \widehat{d}_{Z} ( x , y )
> r .
\end{align*}
Hence $ \widehat{d}_{X} ( x, x^{\prime} ) = \widehat{d}_{Y} ( y,
y^{\prime} ) > r \geq \text{dis} \widehat{ \textit{C}}_{0} . $
This shows that $ \widehat{ \textit{C}}_{0} $ is a $ strong \
correspondence $ between $ X $ and $ Y . $ So $ \inf_{
\widehat{\textit{C }} } ( \text{dis} \widehat{ \textit{C} } ) \leq
\text{dis} \widehat{ \textit{C}}_{0} \leq r . $ Since $ r >
\widehat{d}_{GH} ( X , Y ) $ is arbitrary, we obtain that $ \inf_{
\widehat{\textit{C }} } ( \text{dis} \widehat{ \textit{C} } ) \leq
\widehat{d}_{GH} ( X , Y ) . $
\par \vskip 0.15 cm
Step 2. \ Let $ r > \inf_{ \widehat{\textit{C }} } ( \text{dis}
\widehat{ \textit{C} } ) . $ Then there exists a $ strong \
correspondence $ \ $ \widehat{ \textit{C }}_{0} $ such that $
\text{dis} \widehat{ \textit{C}}_{0} < r . $ Denote $ r_{0} =
\text{dis} \widehat{ \textit{C}}_{0} . $ Let $ Z = X \sqcup Y $ be
the disjoint union, and define a function $ \widehat{d}_{Z} $ on $
Z \times Z $ as follows: \ For $ x \in X $ and $ y \in Y , $
define $$ \widehat{d}_{Z}( x , y ) = \widehat{d}_{Z}( y , x) =
\left\{
\begin{array}{l} r_{0}
\quad \quad \text{if} \ ( x , y ) \in \widehat{ \textit{C}}_{0} ; \\
\widehat{d}_{X}( x , x^{\prime}) \quad \text{if} \ ( x , y )
\notin \widehat{ \textit{C}}_{0} \ \text{and} \ ( x^{\prime} , y )
\in \widehat{ \textit{C}}_{0} \ \text{for some} \ x^{\prime} \in X
.
\end{array}
\right. $$ And $ \widehat{d}_{Z} \mid _{X \times X } =
\widehat{d}_{X} ,  \  \widehat{d}_{Z} \mid _{Y \times Y } =
\widehat{d}_{Y} . $ \\
Firstly, $ \widehat{d}_{Z} $ is well defined. In fact, if $ ( x ,
y ) \notin \widehat{ \textit{C}}_{0} $ with $ x \in X $ and $ y
\in Y , $ then for any $ x^{\prime} , x^{\prime\prime} \in X $
such that $ ( x^{\prime} , y ) ,  ( x^{\prime\prime} , y ) \in
\widehat{ \textit{C}}_{0} , $ by Lemma 2.13.(2), we have $
\widehat{d}_{X}( x , x^{\prime}) = \widehat{d}_{X}( x ,
x^{\prime\prime}) > r_{0} . $ So $ \widehat{d}_{Z} $ is well
defined. Obviously, $ \widehat{d}_{Z} $ is symmetry and $
\widehat{d}_{Z} ( z_{1} , z_{2} ) \geq 0 $ if $ z_{1} \neq z_{2} .
$ In particular, $ \widehat{d}_{Z} ( x , y ) \geq r_{0} $ for any
$ x \in X $ and $ y \in Y . $  \\
Next by the definition and Lemma 2.13.(2), via a direct but
tedious calculation, it can be verified that $ \widehat{d}_{Z} $
satisfies the strong triangle inequality, so $ \widehat{d}_{Z} $
is a non-Archimedean semi-metric on $ Z $ (if $ r_{0} > 0 , $ then
$ \widehat{d}_{Z} $ is a non-Archimedean metric).  As in [BBI,
p.2], we use $ ( Z / \widehat{d}_{Z} , \widehat{d}_{Z} ) $ to
represent the metric space (of course, non-Archimedean) associated
to $ ( Z , \widehat{d}_{Z} ) . $ Then $ X $ and $ Y $ can be
isometrically embedded in $ Z / \widehat{d}_{Z} . $ Now we come to
compute $ \widehat{d}_{Z , H } ( X , Y ) $ in $ Z /
\widehat{d}_{Z} . $ For any $ x \in X , $ by definition, there
exist $ y \in Y $ such that $ ( x , y ) \in \widehat{
\textit{C}}_{0} , $ so $ \widehat{d}_{Z} ( x , y ) = r_{0} , $ and
then $ \text{dist} ( x , Y ) \leq \widehat{d}_{Z} ( x , y ) =
r_{0} < r . $ Hence $ x \in U_{r} ( Y ) $ and so $ X \subset U_{r}
( Y ) . $ Similarly, $ Y \subset U_{r} ( X ) . $ Therefore by
definition $ \widehat{d}_{Z , H } ( X , Y ) \leq r , $ so $
\widehat{d}_{GH } ( X , Y ) \leq \widehat{d}_{Z , H } ( X , Y )
\leq r . $ Since $ r > \inf_{ \widehat{\textit{C }} } ( \text{dis}
\widehat{ \textit{C} } ) $ is arbitrary, we obtain that $
\widehat{d}_{GH } ( X , Y ) \leq \inf_{ \widehat{\textit{C }} } (
\text{dis} \widehat{ \textit{C} } ) . $ Therefore $
\widehat{d}_{GH } ( X , Y ) = \inf_{ \widehat{\textit{C }} } (
\text{dis} \widehat{ \textit{C} } ) . $ This completes the proof
of Theorem 2.14. \quad $ \Box $
\par  \vskip 0.2 cm

{\bf Remark 2.15.} \ From the step 2 of the proof of Theorem 2.14
above, we know that $ \widehat{d}_{Z} $ is a non-Archimedean
metric on $ Z $ if $ r_{0} = \text{dis} \widehat{ \textit{C}}_{0}
> 0 ; $ and $ \widehat{d}_{Z} $ is a semi-metric if $ r_{0} = 0 .
$ The later case is trivial. In fact, if $ r_{0} = 0 , $ then it
is easy to see that $ X $ and $ Y $ are isometric. More precisely,
the map $ f : X \longrightarrow Y $ which sends $ x \in X $ to an
element $ y \in Y $ with $ ( x , y ) \in \widehat{ \textit{C}}_{0}
$ is an isometry, so $ \widehat{d}_{GH } ( X , Y ) = 0. $ From
this we know that, $ (X, \widehat{d}_{X}) $ and $ (Y,
\widehat{d}_{Y}) $ are isometric if and only if there exists a $
strong \ correspondence \ \widehat{ \textit{C}} $ between them
with $ \text{dis} \widehat{ \textit{C}} =  0.  \quad \Box $
\par  \vskip 0.2 cm

{\bf Definition 2.16.} \ For any two non-Archimedean metric spaces
$ ( X , \widehat{d}_{X} ) $ and $ ( Y , \widehat{d}_{Y} ) $ with $
d_{GH } ( X , Y ) > 0 , $ equivalently, $ X $ and $ Y $ are not
isometric. We define a function $$ \Upsilon _{m} ( X ,  Y ) =
\frac{ \widehat{d}_{GH } ( X , Y ) }{ d_{GH } ( X , Y ) } . $$ And
call it the metric ratio function.
\par  \vskip 0.2 cm

Note that by Corollary 2.9, we always have $ \Upsilon _{m} ( X , Y
) \geq 2. $
\par  \vskip 0.2 cm

We will study the properties of the metric ratio function $
\Upsilon _{m} ( X , Y ) $ later. At first we ask the following
\par  \vskip 0.2 cm

{\bf Question 2.17.} \ Is the metric ratio function $ \Upsilon
_{m} $ unbounded ? in other words, for any $ c \geq 2 , $ do there
exist non-Archimedean metric spaces $ X $ and $ Y $ such that $
\widehat{d}_{GH } ( X , Y ) \geq c \cdot d_{GH } ( X , Y ) ? $
\par  \vskip 0.2 cm

We will give an affirmative answer to this question in section 4 (
see Theorem 4.8 below ).
\par  \vskip 0.2 cm

Recall that $ \widehat{\Gamma }_{c} $ is the set of isometry
classes of all compact non-Archimedean metric spaces, and we have
two kinds of metrics $ \widehat{d}_{GH } $ and $ d_{GH } $ on it.
So we ask the following
\par  \vskip 0.2 cm

{\bf Question 2.18.} \ What about the relations between $ (
\widehat{\Gamma }_{c} , \ \widehat{d}_{GH } ) $ and $ (
\widehat{\Gamma }_{c} , \ d_{GH } ) ? $
\par  \vskip 0.2 cm

{\bf Definition 2.19.} \ Let $ \widehat{ \textit{C}} $ be a $
strong \ correspondence $ between two non-Archimedean metric
spaces $ ( X , \widehat{d}_{X} ) $ and $ ( Y , \widehat{d}_{Y} ) .
$ Denote $ \widehat{ \textit{C}}^{\perp } = X \times Y \setminus
\widehat{ \textit{C}} $ be the complement of $ \widehat{
\textit{C}} $ in $ X \times Y . $ For $ \widehat{
\textit{C}}^{\perp } \neq \emptyset , $ we define a function $
\chi _{ \widehat{ \textit{C}} } $ on it as follows: \\
For $ ( x , y ) \in \widehat{ \textit{C}}^{\perp } , $ define $
\chi _{ \widehat{ \textit{C}}} ( x , y ) = \widehat{d}_{Y} ( y ,
y^{\prime} ) $ for some $ y^{\prime} \in Y $ such that $ ( x ,
y^{\prime} ) \in \widehat{ \textit{C}} . $ Then by the following
Lemma 2.20, we know that $ \chi _{ \widehat{ \textit{C}} } $ is
indeed a function on $  \widehat{ \textit{C}}^{\perp } . $ We
simply write it as $ \chi , $ and call it the equilibrium function
associated to $ \widehat{ \textit{C}} $ between $ X $ and $ Y . $
\par  \vskip 0.2 cm

{\bf Lemma 2.20.} \ (1) \ $ \chi _{ \widehat{ \textit{C}} } $ is
well defined and so it is indeed a function on $  \widehat{
\textit{C}}^{\perp } . $
\par  \vskip 0.15 cm
(2) \ For $ ( x , y ) \in \widehat{ \textit{C}}^{\perp } , $ we
have $ \chi _{ \widehat{ \textit{C}}} ( x , y ) = \widehat{d}_{X}
( x , x^{\prime} ) $ for any $ x^{\prime} \in X $ such that $ (
x^{\prime} , y ) \in \widehat{ \textit{C}} . $
\par  \vskip 0.15 cm
(3) \ We have the inequality $$ \text{dis}( \widehat{ \textit{C}}
) < \chi _{ \widehat{ \textit{C}}} ( x , y ) \leq \min \{
\text{diam} ( X ) , \ \text{diam} ( Y ) \} $$ for any $ ( x , y )
\in \widehat{ \textit{C}}^{\perp } . $
\par  \vskip 0.2 cm
{\bf Proof.} \ For any $ y^{\prime\prime} \in Y $ such that $ ( x
, y^{\prime\prime} ) \in \widehat{ \textit{C}} , $ by Lemma
2.13.(2), we have $ \widehat{d}_{Y} ( y , y^{\prime} ) =
\widehat{d}_{Y} ( y , y^{\prime\prime} ) > \text{dis}( \widehat{
\textit{C}} ) . $ So $ \chi _{ \widehat{ \textit{C}}} ( x , y ) =
\widehat{d}_{Y} ( y , y^{\prime} ) = \widehat{d}_{Y} ( y ,
y^{\prime\prime} ) $ is independent of the choice of $ y^{\prime}
. $ Similarly, for any $ ( x^{\prime} , y ) , \ ( x^{\prime\prime}
, y ) \in \widehat{ \textit{C}} $ with $ x^{\prime} ,
x^{\prime\prime} \in X , $ by Lemma 2.13.(2), $ \widehat{d}_{X} (
x , x^{\prime} ) = \widehat{d}_{X} ( x , x^{\prime\prime} ) >
\text{dis}( \widehat{ \textit{C}} ) . $ Moreover, by Lemma 2.12,
we have $ \widehat{d}_{X} ( x , x^{\prime} ) = \widehat{d}_{Y} ( y
, y^{\prime} ) > \text{dis}( \widehat{ \textit{C}} ) . $ On the
other hand, $ \chi _{ \widehat{ \textit{C}}} ( x , y ) =
\widehat{d}_{Y} ( y , y^{\prime} ) =  \widehat{d}_{X} ( x ,
x^{\prime} ) \leq \min \{ \text{diam} ( X ) , \ \text{diam} ( Y )
\}. $ This proves Lemma 2.20. \quad $ \Box $
\par  \vskip 0.2 cm

For this function $ \chi , $ we ask the following
\par  \vskip 0.2 cm

{\bf Question 2.21.} \ (1) \ $$ \text{How about} \inf_{ ( x , \ y
) \in \widehat{ \textit{C}}^{\perp } } \chi_{ \widehat{
\textit{C}}} ( x , y )  \quad \text{and} \quad  \sup_{ ( x , \ y )
\in \widehat{ \textit{C}}^{\perp } } \chi_{ \widehat{ \textit{C}}}
( x , y ) \quad ? $$ Are the former equal to $ \text{dis}(
\widehat{ \textit{C}} ) $ and the later equal to $ \min \{
\text{diam} ( X ) , \ \text{diam} ( Y ) \} \ ? $
\par  \vskip 0.15 cm
(2) \ How do $ \chi_{ \widehat{ \textit{C}}} $ and the quantities
in (1) vary as $ \widehat{ \textit{C}} $ runs over all the $
strong $ \ $ correspondence $ between $ X $ and $ Y  \ ? $
\par  \vskip 0.15 cm
(3) \ If $ X $ and $ Y $ are endowed with some metric measures $
\mu_{X} $ and $ \mu_{Y}, $ then is the function $ \chi_{ \widehat{
\textit{C}}} $ integrable under the product measure ? if so, then
what about the possible integral $ \int_{ \widehat{
\textit{C}}^{\perp } } \chi_{ \widehat{ \textit{C}}} ( x , y ) d
\mu_{X \times Y } \ ? $
\par  \vskip 0.2 cm

We will discuss these questions in a separate paper.
\par  \vskip 0.2 cm

Now we come to establish another formula for explicitly computing
$ \widehat{d}_{GH} . $ To begin with, we define a useful tool, the
$ strong \ \varepsilon -$isometry, as follows:
\par  \vskip 0.2 cm

{\bf Definition 2.22.} \ Let $ ( X , \widehat{d}_{X} ) $ and $ ( Y
, \widehat{d}_{Y} ) $ be two non-Archimedean metric spaces and $
\varepsilon > 0. $ \ Let $ f : X \longrightarrow Y $ be an $
\varepsilon -$isometry. If $ f $ satisfies the following
conditions:
\par \vskip 0.15 cm
( $ SI_{1} $ ) \ For $ x \in X $ and $ y \in Y , $ if $
\widehat{d}_{Y} ( y , f(x) ) \geq \varepsilon , $ then there
exists an $ x^{\prime} \in X $ such that $ \widehat{d}_{Y} ( y ,
f(x^{\prime}) ) < \varepsilon $ and $ \widehat{d}_{X} ( x ,
x^{\prime} ) = \widehat{d}_{Y} ( y , f(x) ) . $
\par \vskip 0.15 cm
( $ SI_{2} $ ) \ For $ x_{1} , \ x_{2} \in X , $ if $
\widehat{d}_{X} ( x_{1} , x_{2} ) \neq \widehat{d}_{Y} ( f( x_{1}
) , f( x_{2}) ), $ then $ \widehat{d}_{X} ( x_{1} , x_{2} ) <
\varepsilon $ \ ( equivalently, if $ \widehat{d}_{X} ( x_{1} ,
x_{2} ) \geq \varepsilon , $ then $ \widehat{d}_{X} ( x_{1} ,
x_{2} ) = \widehat{d}_{Y} ( f( x_{1} ) , f( x_{2}) ) $ ).
\par \vskip 0.15 cm
Then we call $ f $ a $ strong \ \varepsilon -$isometry from $ X $
to $ Y . $
\par  \vskip 0.2 cm

{\bf Theorem 2.23.} \ Let $ ( X , \widehat{d}_{X} ) $ and $ ( Y ,
\widehat{d}_{Y} ) $ be two non-Archimedean metric spaces and $
\varepsilon > 0. $ Then
\par  \vskip 0.15 cm
(1) \ If $ \widehat{d}_{GH} ( X , Y ) < \varepsilon , $ then there
exists a $ strong \ \varepsilon -$isometry from $ X $ to $ Y . $
\par  \vskip 0.15 cm
(2) \ If there exists a $ strong \ \varepsilon -$isometry from $ X
$ to $ Y , $ then $ \widehat{d}_{GH} ( X , Y ) \leq \varepsilon .
$
\par  \vskip 0.2 cm

{\bf Proof.} \ (1) \ If $ \widehat{d}_{GH} ( X , Y ) < \varepsilon
, $ then by Theorem 2.14, \ $ \inf_{\widehat{ \textit{C}}}
(\text{dis} \widehat{ \textit{C}} \ ) < \varepsilon . $ So there
exists a $ strong \ correspondence \ \widehat{ \textit{C}}_{0} $
between $ X $ and $ Y $ such that $ \text{dis} \widehat{
\textit{C}}_{0} < \varepsilon . $ We define a map $ f : X
\longrightarrow Y $ as follows: for each $ x \in X , $ we select
only one $ y \in Y $ satisfying $ ( x , y ) \in \widehat{
\textit{C}}_{0} , $ and then define $ f(x) = y $ ( such $ y $
exists because $ \widehat{ \textit{C}}_{0} $ is a $ correspondence
$ ). Certainly, such $ f $ may not be unique, we fix one among
them. By definition, $ ( x , f(x) ) \in \widehat{ \textit{C}}_{0}
$ for all $ x \in X . $ \\
Step 1. \ By definition
\begin{align*}
\text{dis} f &= \sup_{ x_{1} , \ x_{2} \in X } \mid
\widehat{d}_{Y} ( f(x_{1}) , f(x_{2}) ) - \widehat{d}_{X} ( x_{1}
, x_{2} ) \mid  \\
&\leq \sup \{ \mid \widehat{d}_{Y} ( y_{1} , y_{2} ) -
\widehat{d}_{X} ( x_{1} , x_{2} ) \mid \ : \ (x_{1} , y_{1} ) , \
( x_{2} , y_{2} ) \in \widehat{ \textit{C}}_{0} \} \\
&= \text{dis} \widehat{ \textit{C}}_{0} < \varepsilon . \quad
\quad  \text{That is,} \  \text{dis} f < \varepsilon .
\end{align*}
Step 2. \ For any $ y \in Y , $ there exists $ x \in X $ such that
$ ( x , y ) \in \widehat{ \textit{C}}_{0}. $ Because $ ( x , f(x)
) \in \widehat{ \textit{C}}_{0}, $ we have $$ \widehat{d}_{Y} ( y
, f(x) ) = \mid \widehat{d}_{Y} ( y , f(x) ) - \widehat{d}_{X} ( x
, x ) \mid \ \leq \text{dis} \widehat{ \textit{C}}_{0} <
\varepsilon . $$ So $ \text{dist} ( y , f ( X )) \leq
\widehat{d}_{Y} ( y , f(x) ) < \varepsilon . $ Since $ y $ is
arbitrary, this shows that
$ f ( X ) $ is an $ \varepsilon -$net in $ Y . $ \\
The above two steps show that $ f $ is an $ \varepsilon -$isometry
from $ X $ to $ Y . $ \\
Step 3. \ Let $ x \in X $ and $ y \in Y. $ If $ \widehat{d}_{Y} (
y , f(x) ) \geq \varepsilon , $ then $ ( x , y ) \notin \widehat{
\textit{C}}_{0}. $ Otherwise, the fact that both $ ( x , f(x)) \in
\widehat{ \textit{C}}_{0} $ and $ ( x , y ) \in \widehat{
\textit{C}}_{0} $ imply $$ \widehat{d}_{Y} ( y , f(x) ) = \mid
\widehat{d}_{Y} ( y , f(x) ) - \widehat{d}_{X} ( x , x ) \mid \
\leq \text{dis} \widehat{ \textit{C}}_{0} < \varepsilon . $$ A
contradiction! So $ ( x , y ) \notin \widehat{ \textit{C}}_{0}. $
Then by Lemma 2.12, for any $ x^{\prime} \in X $ satisfying $ (
x^{\prime} , y ) \in \widehat{ \textit{C}}_{0}, $ we have $
\widehat{d}_{X} ( x , x^{\prime} ) = \widehat{d}_{Y} ( y , f(x) )
. $ \ Since $ ( x^{\prime} , f(x^{\prime}) ) \in \widehat{
\textit{C}}_{0} , $ we have
$$ \widehat{d}_{Y} ( y , f(x^{\prime}) ) = \mid \widehat{d}_{Y}
( y , f(x^{\prime}) ) - \widehat{d}_{X} ( x^{\prime} , x^{\prime}
) \mid \ \leq \text{dis} \widehat{ \textit{C}}_{0} < \varepsilon .
$$ Since such $ x^{\prime} $ always exists, this shows that $ f $
satisfies the condition $ ( SI_{1}) $ of Def.2.22. \\
Step 4. \ Let $ x_{1} , x_{2} \in X . $ If $ \widehat{d}_{X} (
x_{1} , x_{2} ) \neq \widehat{d}_{Y} ( f(x_{1}) , f(x_{2}) ), $
then $ ( x_{2} , f(x_{1}) ) \in \widehat{ \textit{C}}_{0} . $
Otherwise, by the fact that $ ( x_{2} , f(x_{1}) ) \notin
\widehat{ \textit{C}}_{0} $ and $ ( x_{i} , f(x_{i}) ) \in
\widehat{ \textit{C}}_{0} \ ( i = 1 , 2 ), $ we get via Lemma 2.12
that $ \widehat{d}_{X} ( x_{1} , x_{2} ) = \widehat{d}_{Y} (
f(x_{1}) , f(x_{2}) ), $ a contradiction! so $ ( x_{2} , f(x_{1})
) \in \widehat{ \textit{C}}_{0}. $ Since $  ( x_{1} , f(x_{1}) )
\in \widehat{ \textit{C}}_{0}, $ we have $$ \widehat{d}_{X} (
x_{1} , x_{2} ) = \mid \widehat{d}_{X} ( x_{1} , x_{2} ) -
\widehat{d}_{Y} ( f(x_{1}) , f(x_{1}) ) \mid \ \leq \text{dis}
\widehat{ \textit{C}}_{0} < \varepsilon . $$ This shows that $ f $
satisfies the condition $ ( SI_{2}) $ of Def.2.22. \\
To sum up, $ f $ is a $ strong \ \varepsilon -$isometry from $ X $
to $ Y. $ This proves (1).
\par  \vskip 0.15 cm
(2) \ Let $ f : X \longrightarrow Y $ be a $ strong \ \varepsilon
-$isometry. Define a subset $ \widehat{ \textit{C}}_{0} ( f )
\subset X \times Y $ by $$ \widehat{ \textit{C}}_{0} ( f ) = \{ (
x , y ) \in X \times Y  \ : \ \widehat{d}_{Y} ( y, f(x) ) \leq
\varepsilon \}. $$ Step 1$^{\prime}. $ \ For any $ x \in X , $
obviously $ ( x , f(x) ) \in \widehat{ \textit{C}}_{0} ( f ). $
For any $ y \in Y, $ since $ f(X) $ is an $ \varepsilon -$net in $
Y , \ \text{dist}( y, f(X)) < \varepsilon , $ so there exists $ x
\in X $ such that $ \widehat{d}_{Y} ( y, f(x) ) < \varepsilon , $
and then $ ( x , y ) \in \widehat{ \textit{C}}_{0} ( f ). $ So $
\widehat{ \textit{C}}_{0} ( f ) $ is
a $ correspondence $ between $ X $ and $ Y . $  \\
Step 2$^{\prime}. $ \ Let $ ( x_{1} , y_{1} ) , \ ( x_{2} , y_{2}
) \in \widehat{ \textit{C}}_{0} ( f ). $ Then $ \widehat{d}_{Y} (
y_{1}, f(x_{1}) ) \leq \varepsilon  $ and $ \widehat{d}_{Y} (
y_{2}, f(x_{2}) ) \leq \varepsilon . $ \\
Case A. \ We assume that $ \widehat{d}_{Y} ( f(x_{1}), f(x_{2}) )
> \varepsilon . $ Then by the strong triangle inequality, we have
\begin{align*}
\widehat{d}_{Y} ( y_{1} , y_{2} ) &= \max \{ \widehat{d}_{Y} (
y_{1} , f(x_{1}) ) , \ \widehat{d}_{Y} ( f(x_{1}), f(x_{2}) ), \
\widehat{d}_{Y} ( f(x_{2}), y_{2} ) \} \\
&= \widehat{d}_{Y} ( f(x_{1}), f(x_{2}) ) .  \quad \quad \text{ So
}
\end{align*}
$$ \mid \widehat{d}_{X} ( x_{1} , x_{2} ) - \widehat{d}_{Y} ( y_{1} ,
y_{2} ) \mid \ = \ \mid \widehat{d}_{X} ( x_{1} , x_{2} ) -
\widehat{d}_{Y} ( f(x_{1}), f(x_{2}) ) \mid \ \leq \text{dis} f <
\varepsilon . $$ Case B. \ We assume that $ \widehat{d}_{Y} (
f(x_{1}), f(x_{2}) ) \leq \varepsilon , $ and discuss this case
via the following two subcases: \\
Subcase B1. \ We assume that $ \widehat{d}_{X} ( x_{1}, x_{2} ) =
\widehat{d}_{Y} ( f(x_{1}), f(x_{2})). $ Then since $$
\widehat{d}_{Y} ( y_{1} , y_{2} ) \leq \max \{ \widehat{d}_{Y} (
y_{1} , f(x_{1}) ) , \ \widehat{d}_{Y} ( f(x_{1}), f(x_{2}) ), \
\widehat{d}_{Y} ( f(x_{2}), y_{2} ) \} \leq \varepsilon , $$ we
have $ \mid \widehat{d}_{Y} ( y_{1} , y_{2} ) - \widehat{d}_{Y} (
f(x_{1}), f(x_{2}) ) \mid \ \leq \varepsilon , $ so $ \mid
\widehat{d}_{X} ( x_{1} , x_{2} ) - \widehat{d}_{Y} ( y_{1}, y_{2}
) \mid \ \leq \varepsilon . $ \\
Subcase B2. \ We assume that $ \widehat{d}_{X} ( x_{1}, x_{2} )
\neq \widehat{d}_{Y} ( f(x_{1}), f(x_{2})). $ Then by the fact
that $ f $ is a $ strong \ \varepsilon -$isometry we get $
\widehat{d}_{X} ( x_{1}, x_{2} ) < \varepsilon $ ( via the
condition ( $ SI_{2} $ ) ). As above, $ \widehat{d}_{Y} ( y_{1} ,
y_{2} ) \leq \varepsilon . $ So $ \mid \widehat{d}_{X} ( x_{1} ,
x_{2} ) - \widehat{d}_{Y} ( y_{1}, y_{2} ) \mid \ \leq \varepsilon
. $  \\
To sum up, we always have $ \mid \widehat{d}_{X} ( x_{1} , x_{2} )
- \widehat{d}_{Y} ( y_{1}, y_{2} ) \mid \ \leq \varepsilon  $ for
all $ ( x_{1} , y_{1} ), \ ( x_{2} , y_{2} ) \in \widehat{
\textit{C}}_{0} ( f ). $ This shows $ \text{dis} \widehat{
\textit{C}}_{0} ( f ) \leq \varepsilon . $ \\
Step 3$^{\prime}$. \ For $ x \in X $ and $ y \in Y , $ if $ ( x ,
y ) \notin \widehat{ \textit{C}}_{0} ( f ), $ then $
\widehat{d}_{Y} ( y , f(x)) > \varepsilon . $ So by the condition
( $ SI_{1} $ ) for $ f , $ there exists an $ x^{\prime} \in X $
such that $ \widehat{d}_{Y} ( y , f(x^{\prime})) < \varepsilon $
and $ \widehat{d}_{X} ( x , x^{\prime} ) = \widehat{d}_{Y} ( y ,
f(x)) . $ In particular,  $ ( x^{\prime} , y ) \in \widehat{
\textit{C}}_{0} ( f ) . $ Moreover, $ ( x , f(x)) \in \widehat{
\textit{C}}_{0} ( f ) . $ This shows $ \widehat{ \textit{C}}_{0} (
f ) $ satisfying the condition ( $ C_{NA} $ ) of Def. 2.11, hence
$ \widehat{ \textit{C}}_{0} ( f ) $ is a $ strong \ correspondence
$ between $ X $ and $ Y . $ Therefore by Theorem 2.14, we obtain $
\widehat{d}_{GH} ( X , Y ) \leq \text{dis} \widehat{
\textit{C}}_{0} ( f ) \leq \varepsilon . $ This proves (2), and
the proof of Theorem 2.23 is completed. \quad $ \Box $
\par  \vskip 0.2 cm

{\bf Corollary 2.24.} \ Let $ ( X , \widehat{d}_{X} ) $ and $ ( Y
, \widehat{d}_{Y} ) $ be two non-Archimedean metric spaces. Then
$$ \widehat{d}_{GH} ( X , Y ) = \inf \{ \varepsilon > 0 \ : \
\text{there exists a } \ strong \ \varepsilon -\text{isometry
from} \ X \ \text{to} \ Y  \} . $$
\par  \vskip 0.1 cm
{\bf Proof.} \ Denote $ d = \inf \{ \varepsilon > 0 \ : \
\text{there exists a } \ strong \ \varepsilon -\text{isometry
from} \ X \ \text{to} \ Y  \} . $ For any $ \varepsilon >
\widehat{d}_{GH} ( X , Y ) , $ by Theorem 2.23.(1), there exists a
$ strong $ \ $ \varepsilon -$isometry from $ X $ to $ Y, $ so $ d
\leq \varepsilon , $ and then $ d \leq \widehat{d}_{GH} ( X , Y )
. $ \\
Conversely, for any $ \varepsilon > d , $ there exists a $ strong
$ \ $ \varepsilon ^{\prime}-$isometry from $ X $ to $ Y $ with $ d
\leq \varepsilon ^{\prime} < \varepsilon . $ So by Theorem
2.23.(2), we have $ \widehat{d}_{GH} ( X , Y ) \leq \varepsilon
^{\prime} < \varepsilon . $ Since $ \varepsilon > d $ is
arbitrary, we get $ \widehat{d}_{GH} ( X , Y ) \leq d . $ Hence $
\widehat{d}_{GH} ( X , Y ) = d , $ and the proof is completed.
\quad $ \Box $

\par     \vskip  1 cm

\hspace{-0.6cm}{\bf 3. \ Strong Gromov-Hausdorff convergence }

\par \vskip 0.8 cm

In this section, we consider the converging sequences in the
collection of non-Archimedean metric spaces, especially in the
non-Archimedean Gromov-Hausdorff space $ ( \widehat{\Gamma }_{c} ,
\ \widehat{d}_{GH } ). $ Because there are two metrics $
\widehat{d}_{GH } $ and $ d_{GH } $ as defined, to distinguish the
convergence corresponding to them, we call the convergence under $
\widehat{d}_{GH } $ the strong type, as stated precisely in the
following
\par  \vskip 0.2 cm

{\bf Definition 3.1.} \ A sequence $ \{ X_{n} \}_{n \in \N } $ of
non-Archimedean metric spaces  strongly converges to a
non-Archimedean metric space $ X $ if $ \widehat{d}_{GH } ( X_{n}
,  X) \longrightarrow 0 $ as $ n \longrightarrow \infty . $ In
this case, we will write $ X_{n} \longrightarrow_{ \text{GH}_{S}}
X $ and call $ X $ a strong  Gromov-Hausdorff limit of $ \{ X_{n}
\} . $ If all $ X_{n} $ and $ X $ are compact, then the limit is
unique up to an isometry since $ \widehat{d}_{GH } $ is a metric
on $ \widehat{\Gamma }_{c} . $
\par  \vskip 0.15 cm
As usual, if $ d_{GH } ( X_{n} ,  X) \longrightarrow 0 $ as $ n
\longrightarrow \infty , $ then we write $ X_{n} \longrightarrow_{
\text{GH}} X $ and call $ X $ a  Gromov-Hausdorff limit of $ \{
X_{n} \}  $ ( see [BBI, p.260] ).
\par  \vskip 0.15 cm
If $ X_{n} \longrightarrow_{ \text{GH}_{S}} X , $ then by
Corollary 2.9, we have $ X_{n} \longrightarrow_{ \text{GH}} X . $
\par  \vskip 0.2 cm

{\bf Example 3.2.} \ Every compact non-Archimedean metric space $
( X , \widehat{d}_{X}) $ is a strong Gromov-Hausdorff limit of
finite spaces.
\par  \vskip 0.2 cm

{\bf Proof.} \ This is a counterpart of Example 7.4.9 of [BBI],
the proof is similar. \quad $ \Box $ \\
If we write $ \widehat{\Gamma }_{c , f } = \{ X : \ X \in
\widehat{\Gamma }_{c} \text{ and } X \text{ is a finite set } \},
$ then Example 3.2 shows that $ \widehat{\Gamma }_{c , f } $ is
dense in $ \widehat{\Gamma }_{c} $ under the metric $
\widehat{d}_{GH}. $
\par  \vskip 0.15 cm
Likewise, as a corollary of Theorem 2.23, there is a similar
criterion for the strong convergence corresponding to the usual
one ( see [BBI, p.260]) as follows:
\par  \vskip 0.2 cm

{\bf Criterion 3.3.} \ A sequence $ \{ X_{n} \}_{n \in \N } $ of
non-Archimedean metric spaces  strongly converges to a
non-Archimedean metric space $ X $ if and only if there are a
sequence $ \{ \varepsilon _{n} \} $ of positive real numbers and a
sequence of maps $ f_{n} : X_{n} \longrightarrow X $ ( or
alternatively, $ f_{n} : X \longrightarrow X_{n} $ ) such that
every $ f_{n} $ is a $ strong \ \varepsilon _{n}-$isometry and $
\varepsilon _{n} \longrightarrow 0 . $
\par  \vskip 0.2 cm

{\bf Proof.} \ $ \Longleftarrow : $ \ For any $ \varepsilon > 0 ,
$ there exist $ n_{0} \in \N $ such that $ \varepsilon _{n} <
\varepsilon $ for all $ n > n_{0} . $ Since every $ f_{n} : X_{n}
\longrightarrow X $ is a strong $ \varepsilon _{n} -$isometry, by
Theorem 2.23.(2), $ \widehat{d}_{GH }( X_{n} ,  X) \leq
\varepsilon _{n} < \varepsilon , $ so
$ X_{n} \longrightarrow_{ \text{GH}_{S}} X . $ \\
$ \Longrightarrow : $ \ If $ X_{n} \longrightarrow_{
\text{GH}_{S}} X , $ then for any $ \varepsilon > 0 , $ there
exist $ n_{0} \in \N
 $ such that $ \widehat{d}_{GH } ( X_{n} ,  X) < \varepsilon  $ for all
 $ n >  n_{0}. $ So by Theorem 2.23.(1), there exists a strong
 $ \varepsilon -$isometry $ f_{n} : X \longrightarrow X_{n} $ for
 every $ n >  n_{0} . $ By letting $ \varepsilon \longrightarrow 0
 $ as $ n \longrightarrow \infty , $ we obtain the results. \quad
 $ \Box $
\par  \vskip 0.2 cm

Now we define another tool, the $ strong \ \varepsilon
-approximation , $ which will also be used in explicit computation
of $ \widehat{d}_{GH } . $
\par  \vskip 0.2 cm

{\bf Definition 3.4.} \ Let $ ( X , \widehat{d}_{X} ) $ and $ ( Y
, \widehat{d}_{Y} ) $ be two compact non-Archimedean metric spaces
and $ \varepsilon > 0 . $ We say $ X $ and $ Y $ are $ strong \
\varepsilon -approximations $ of each other if there exist finite
collections of points $ \{ x_{i} \}_{i =  1 }^{N} $ and $ \{ y_{i}
\}_{i =  1 }^{N} $ in $ X $
 and $ Y , $ respectively, such that :
\par  \vskip 0.15 cm
(1) \ The set $ \{ x_{i} \ : \ 1\leq i \leq N \} $ is an $
\varepsilon -$net in $ X , $ and $ \{ y_{i} \ : \ 1\leq i \leq N
\} $ is an $ \varepsilon -$net in $ Y. $
\par  \vskip 0.15 cm
(2) \ $ \widehat{d}_{X} ( x_{i} , \ x_{j} ) =  \widehat{d}_{Y} (
y_{i} , \ y_{j} ) $ for all $ i , \ j \in \{ 1, \cdots , N \} . $
\par  \vskip 0.2 cm
Obviously, if $ X $ and $ Y $ are $ strong \ \varepsilon
-approximations $ of each other, then they are also $ \varepsilon
-approximations $ of each other in the sense of [BBI,Def.7.4.10].
\par  \vskip 0.2 cm

{\bf Theorem 3.5.} \ Let $ ( X , \widehat{d}_{X} ) $ and $ ( Y ,
\widehat{d}_{Y} ) $ be two compact non-Archimedean metric spaces.
\par  \vskip 0.15 cm
(1) \ If $ X $ and $ Y $ are $ strong \ \varepsilon
-approximations $ of each other, then $ \widehat{d}_{GH} ( X , Y )
\leq \varepsilon . $
\par  \vskip 0.15 cm
(2) \ If $ \widehat{d}_{GH} ( X , Y ) < \varepsilon , $ then $ X $
and $ Y $ are $ strong \ \varepsilon -approximations $ of each
other.
\par  \vskip 0.2 cm

{\bf Proof.} \ (1) \ Since $ X $ and $ Y $ are $ strong \
\varepsilon -approximations $ of each other, by definition, there
exist finite collections of points $ \{ x_{i} \}_{i =  1 }^{N}
\subset  X $ and $ \{ y_{i} \}_{i =  1 }^{N} \subset  Y $ such
that conditions (1) and (2) in Def.3.4 hold. Denote $ X_{0} = \{
x_{i} \}_{i =  1 }^{N} $ and $ Y_{0} = \{ y_{i} \}_{i =  1 }^{N},
$ and define $ \textit{C}_{0} = \{ ( x_{i} , y_{i}) \ : \ 1 \leq i
\leq N  \} \subset X_{0} \times Y_{0} . $ Obviously, $
\textit{C}_{0} $ is a $ correspondence $ between $ X_{0} $ and $
Y_{0}, $ and $ \text{dis} ( \textit{C}_{0} ) = 0 $ via condition
(2) of Def.3.4. Now for any $ ( x , y ) \in X_{0} \times Y_{0}
\setminus \textit{C}_{0}, $ we have $ x = x_{i} $ and $ y = y_{j}
$ for some $ i , j \in \{ 1 , \cdots N \} $ and $ i \neq j , $ so
$ x_{i} \neq x_{j} $ and $ y_{i} \neq y_{j} , $ and then $
\widehat{d}_{X_{0}} ( x_{i} , x_{j} ) = \widehat{d}_{Y_{0}} (
y_{i} , y_{j} ) > 0 = \text{dis} ( \textit{C}_{0} ) . $ This shows
that $ \textit{C}_{0} $ is a $ strong \ correspondence $ between $
X_{0} $ and $ Y_{0}. $ Hence by Theorem 2.14, $ \widehat{d}_{GH} (
X_{0} , Y_{0} ) \leq \text{dis} ( \textit{C}_{0} ) = 0 . $ So $
\widehat{d}_{GH} ( X_{0} , Y_{0} ) = 0 . $ Since $ X_{0} $ and $
Y_{0} $ are $ \varepsilon -$nets in $ X $ and $ Y , $
respectively, we have $ \text{dist} ( x , X_{0} ) < \varepsilon  $
and $ \text{dist} ( y , Y_{0} ) < \varepsilon  $ for all $ x \in X
$ and $ y \in Y . $  So
\begin{align*}
&\widehat{d}_{X , H} ( X_{0} , X ) = \max \{ \sup_{ x \in X}
\text{dist} ( x , X_{0} ) , \ \sup_{ x \in X_{0}} \text{dist} ( x
, X ) \} \leq \varepsilon , \quad \text{and} \\
&\widehat{d}_{Y , H} ( Y_{0} , Y ) = \max \{ \sup_{ y \in Y}
\text{dist} ( y , Y_{0} ) , \ \sup_{ y \in Y_{0}} \text{dist} ( y,
Y ) \} \leq \varepsilon .
\end{align*}
Hence $ \widehat{d}_{GH} ( X_{0} , X ) \leq \widehat{d}_{X , H} (
X_{0} , X ) \leq \varepsilon $ \ and \ $ \widehat{d}_{GH} ( Y_{0}
, Y) \leq \widehat{d}_{Y , H} ( Y_{0} , Y ) \leq \varepsilon . $
Therefore by prop.2.6 ( Zarichnyi), we get $$ \widehat{d}_{GH} ( X
, Y ) \leq \max \{ \widehat{d}_{GH} ( X , X_{0} ) , \
\widehat{d}_{GH} ( X_{0} , Y_{0} ) , \ \widehat{d}_{GH} ( Y_{0} ,
Y ) \} \leq \varepsilon . $$ This proves (1).
\par  \vskip 0.15 cm
(2) \ If $ \widehat{d}_{GH} ( X , Y ) < \varepsilon , $ then by
Theorem 2.23.(1), there exists a $ strong \ \varepsilon -$isometry
$ f : X \longrightarrow Y . $ For any $ x \in X, $ denote the ball
$ B_{x} (\varepsilon ) = \{ x^{\prime} \in X \ : \ \widehat{d}_{X}
(x , x^{\prime} ) < \varepsilon  \} , $ which is both open and
closed in $ X $ ( see [Sc]). Then $ X = \cup_{x \in X} B_{x}
(\varepsilon ).$ Since $ X $ is compact and non-Archimedean, we
have $ X = \sqcup_{i = 1 }^{N} B_{x_{i}} (\varepsilon ) $ ( the
disjoint union ) for a positive integer $ N . $ Let $ X_{0} = \{
x_{i} : \ 1 \leq i \leq N  \} $  and define $ y_{i} = f(x_{i}) \in
Y $ for every $ i \in \{  1,  \cdots ,  N  \}. $ We denote $ Y_{0}
= \{ y_{i} : \ 1 \leq i \leq N  \}. $ \\
Now for any $ x_{i}, \ x_{j} \in X_{0} $ with $ i \neq j , $ since
$ B_{x_{i}} (\varepsilon ) \cap B_{x_{j}} (\varepsilon ) =
\emptyset , $ we have $ \widehat{d}_{X} (  x_{i}, \ x_{j} ) \geq
\varepsilon . $ Then by the fact that $ f $ is a strong $
\varepsilon -$isometry we get $ \widehat{d}_{X} (  x_{i}, \ x_{j}
) = \widehat{d}_{Y} (  f(x_{i}), \ f(x_{j}) ) = \widehat{d}_{Y} (
y_{i}, \ y_{j} ). $ So the condition (2) of Def.3.4 holds for $
X_{0} $ and $ Y_{0}. $ \\
Moreover, for any $ x \in X, $ we have $ x \in B_{x_{i}}
(\varepsilon ) $ for some $ i \in \{ 1 , \cdots , N \}, $ so $
\widehat{d}_{X} (  x, \ x_{i} ) < \varepsilon . $ Then $
\text{dist} ( x , X_{0}) \leq \widehat{d}_{X} (  x, \ x_{i} ) <
\varepsilon . $ This shows that $ X_{0} $ is an $ \varepsilon
-$net in $ X . $ \\
Furthermore, since $ f(X) $ is an $ \varepsilon -$net in $ Y , $
for any $ y \in Y, $ there exists an $ x \in X $ such that $
\widehat{d}_{Y} ( y , f(x) ) < \varepsilon . $ Let $ x \in
B_{x_{i}} (\varepsilon ), $ then $ \widehat{d}_{X} ( x , x_{i} ) <
\varepsilon . $  If $ \widehat{d}_{Y} ( f(x) , f(x_{i}) ) \geq
\varepsilon , $ then by the condition $ (SI_{1}) $ of Def.2.22,
there exists an $ x^{\prime} \in X $ such that $ \widehat{d}_{Y} (
f(x) , f(x^{\prime})) < \varepsilon $ and $ \widehat{d}_{X} (
x_{i} , x^{\prime}) = \widehat{d}_{Y} ( f(x) , f(x_{i})) ( \geq
\varepsilon ). $ But $ \widehat{d}_{X} ( x, x^{\prime}) = \max \{
\widehat{d}_{X} ( x, x_{i}) , \widehat{d}_{X} ( x_{i} ,
x^{\prime}) \} = \widehat{d}_{X} ( x_{i} , x^{\prime}) \geq
\varepsilon , $ so by the condition $ (SI_{2}) $ of Def.2.22, we
get $ \widehat{d}_{Y} ( f(x) , f(x^{\prime})) = \widehat{d}_{X} (
x, x^{\prime}) \geq \varepsilon , $ a contradiction! So we must
have $ \widehat{d}_{Y} ( f(x) , f(x_{i})) < \varepsilon . $ Then
$$ \widehat{d}_{Y} (y,  y_{i}) = \widehat{d}_{Y} (y, f(x_{i}))
\leq \max \{ \widehat{d}_{Y} ( y, f(x)) , \ \widehat{d}_{Y} (
f(x), f(x_{i})) \} < \varepsilon . $$ So $ \text{dist} ( y, Y_{0})
 \leq \widehat{d}_{Y} (y,  y_{i}) < \varepsilon . $ This shows
 that $ Y_{0} $ is an $ \varepsilon -$net in $ Y. $  So the condition
 (1) of Def.3.4 holds for $ X_{0}, Y_{0}. $  Therefore, $ X $ and $ Y
$ are $ strong \ \varepsilon - approximations $ of each other.
This proves (2), and the proof of Theorem 3.5 is completed. \quad
$ \Box $
\par  \vskip 0.2 cm

{\bf Corollary 3.6.} \ For compact non-Archimedean metric spaces $
\{ X_{n} \}_{n = 1}^{\infty } $ and $ X , \ X_{n} \longrightarrow
_{\text{GH}_{S}} X $ if and only if, for any $ \varepsilon > 0 , \
X_{n} $ and $ X $ are $ strong \ \varepsilon - approximations $ of
each other for all sufficiently large $ n. $
\par  \vskip 0.2 cm

{\bf Proof.} \ If $  X_{n} \longrightarrow _{\text{GH}_{S}} X , $
then for any $ \varepsilon > 0 , $ there exists a $ n_{0} \in \N $
such that $ \widehat{d}_{GH} ( X_{n} , X ) < \varepsilon $ for all
$ n > n_{0} . $ So by Theorem 3.5.(2), $ X_{n} $ and $ X $ are $
strong \ \varepsilon - approximations $ of each other for all
$ n > n_{0} . $ \\
Conversely, if for any $ \varepsilon > 0 , \ X_{n} $ and $ X $ are
$ strong \ \varepsilon - approximations $ of each other for all
sufficiently large $ n, $ then by Theorem 3.5.(1), \ $
\widehat{d}_{GH} ( X_{n} , X ) \leq \varepsilon $ for all
sufficiently large $ n. $ This shows that $  X_{n} \longrightarrow
_{\text{GH}_{S}} X . $ \quad $ \Box $
\par  \vskip 0.2 cm

{\bf Corollary 3.7.} \ For two compact non-Archimedean metric
spaces $ ( X , \widehat{d}_{X} ) $ and $ ( Y , \widehat{d}_{Y} ),
$ \ $$ \widehat{d}_{GH} ( X , Y ) = \inf \{ \varepsilon > 0 \ : \
X \ \text{and} \ Y \ \text{are} \ strong \ \varepsilon
-aprroximations \ \text{of each other } \}. $$
\par  \vskip 0.2 cm
{\bf Proof.} \ Denote $$ d = \inf \{ \varepsilon > 0 \ : \ X \
\text{and} \ Y \ \text{are} \ strong \ \varepsilon -aprroximations
\ \text{of each other } \}. $$ For any $ \varepsilon >
\widehat{d}_{GH} ( X , Y ) , $ by Theorem 3.5.(2), $ X $ and $ Y $
are $ strong \ \varepsilon - approximations $ of each other, so $
d \leq \varepsilon , $ and then $ d \leq \widehat{d}_{GH} ( X
, Y ). $ \\
Conversely, for any $ \varepsilon > d , $ there exists $ 0 <
\varepsilon ^{\prime} < \varepsilon $ such that $ X $ and $ Y $
are $ strong \ \varepsilon ^{\prime}- approximations $ of each
other. So by Theorem 3.5.(1), we have $ \widehat{d}_{GH} ( X , Y )
 \leq \varepsilon ^{\prime} < \varepsilon , $ this implies that
$ \widehat{d}_{GH} ( X , Y ) \leq d. $ Hence $ \widehat{d}_{GH} (
X , Y ) = d . $ \quad $ \Box $
\par  \vskip 0.2 cm

Based on Theorem 3.5, we can now establish the following Theorem
of determining $ X_{n} \longrightarrow _{\text{GH}_{S}} X , $
which is well compatible to Proposition 7.4.12 of [BBI] about the
usual $ X_{n} \longrightarrow _{\text{GH}} X. $
\par  \vskip 0.2 cm

{\bf Theorem 3.8.} \ For compact non-Archimedean metric spaces $ X
$ and $ \{ X_{n} \}_{n = 1}^{\infty }, $  \ $  X_{n}
\longrightarrow _{\text{GH}_{S}} X $ if and only if the following
holds: \ for every $ \varepsilon > 0 , $ there exists a finite $
\varepsilon -$net $ S $ in $ X $ and a finite $ \varepsilon -$net
$ S_{n} $ in each $ X_{n} $ such that $ S_{n} \longrightarrow
_{\text{GH}_{S}} S . $ Moreover, these $ \varepsilon -$nets can be
chosen so that, for all sufficiently large $ n, \ S_{n} $ have the
same cardinality as $ S. $
\par  \vskip 0.2 cm

Note that for compact non-Archimedean metric spaces, Theorem 3.8
here and Proposition 7.4.12 of [BBI] are not equivalent, this is
because $ \widehat{d}_{GH} $ and $ d_{GH} $ are not equivalent
(see Theorem 4.8 below). Difference between them can also be seen
in their proofs.
\par  \vskip 0.2 cm

{\bf Proof of Theorem 3.8.} \ $ \Longleftarrow : $ \ If such $
\varepsilon -$nets exist, then by Corollary 3.6,  the fact that $
S_{n} \longrightarrow _{\text{GH}_{S}} S  $ implies, for any $
\delta > 0 , $ that $ S_{n} $ and $ S $ are $ strong \ \delta -
approximations $ of each other for all sufficiently large $ n. $
Denote
\begin{align*}
&r_{n} = \min \{ \widehat{d}_{n} ( x , \ x^{\prime}) \ : \ x , \
x^{\prime} \in S_{n} \ \text{and} \ x \neq  x^{\prime} \}, \\
&r = \min \{ \widehat{d}( x , \ x^{\prime}) \ : \ x , \ x^{\prime}
\in S \ \text{and} \ x \neq  x^{\prime} \},
\end{align*}
where $ \widehat{d}_{n} $ and $ \widehat{d} $ are the metrics of $
X_{n} $ and $ X , $ respectively. Then both $ r_{n} > 0 $ and $ r
> 0 $ because $ S_{n} $ and $ S $ are finite sets. Now we take a
$ \delta $ such that $ 0 < \delta <  \min \{ r , r_{n} \}. $ Since
$ S_{n} $ and $ S $ are $ strong \ \delta -approximations $ of
each other for all sufficiently large $ n, $ there exist finite
collections of points $ \{x_{i}\}_{i = 1 }^{N} $ and $
\{y_{i}\}_{i = 1 }^{N} $ in $ S_{n} $ and $ S , $ respectively,
such that: \\
(1) \ The set $ X_{0} = \{ x_{i} \ : \ 1 \leq i \leq N \} $ is a $
\delta -$net in $ S_{n} , $ and $ Y_{0} = \{ y_{i} \ :
\ 1 \leq i \leq N \} $ is a $ \delta -$net in $ S . $ \\
(2) \ $ \widehat{d}_{n} ( x_{i}, \ x_{j}) = \widehat{d} ( y_{i}, \
y_{j}) $ \ for all $ i, \ j \in \{ 1 , \cdots , N \} . $ \\
Then for every $ x \in S_{n}, \ \text{dist} ( x , X_{0}) < \delta
 . $ So there exists an $ x_{i} \in X_{0} $ such that
 $ \widehat{d}_{n} ( x, \ x_{i}) < \delta < r_{n} . $ Thus $ x = x_{i}
 $ because we would have $ r_{n} \leq \widehat{d}_{n} ( x, \ x_{i})
 $ if $ x \neq x_{i} . $ Hence $ S_{n} = X_{0} . $ \\
 Similarly, for every $ y \in S, \ \text{dist} ( y , Y_{0}) <
 \delta , $ \ so $ \widehat{d} ( y, \ y_{i}) < \delta < r $ for
 some $ y_{i} \in Y_{0} . $ Thus $ y = y_{i} \in Y_{0} , $ so
 $ S = Y_{0}. $ These facts together with the above condition (2)
 shows that $ \sharp S_{n} = \sharp X_{0} = \sharp Y_{0} = \sharp S
 $ for all sufficiently large $ n. $ \\
Moreover, since $ S_{n} $ and $ S $ are finite $ \varepsilon
-$nets in $ X_{n} $ and $ X, $ respectively, the above discussion
shows that $ X_{n} $ and $ X $ are $ strong \ \varepsilon
-approximations $ of each other for all sufficiently large $ n. $
Hence by Corollary 3.6, $  X_{n} \longrightarrow _{\text{GH}_{S}}
X . $
\par  \vskip 0.15 cm
$ \Longrightarrow : $ \ If $  X_{n} \longrightarrow
_{\text{GH}_{S}} X , $ then for any $ \varepsilon > 0 , \
\widehat{d}_{GH}( X_{n} , X ) < \varepsilon $ for all sufficiently
large $ n. $ Take a finite $ \varepsilon -$net $ S $ in $ X $ (
such $ S $ exists since $ X $ is compact, hence totally bounded ).
We may as well write $ S = \{ x_{i} \ : \ 1 \leq i \leq N  \} $
with $ N = \sharp S < \infty . $ Denote $$ r = \min \{ \widehat{d}
( x , \ x^{\prime} ) \ : \ x, \ x^{\prime} \in S \ \text{and} \ x
\neq x^{\prime} \}. $$ Then $ r > 0 $ since $ S $ is a finite set.
By Criterion 3.3, there are a sequence of $ strong \ \varepsilon
_{n}-$isometries $ f_{n} \ : \ X \longrightarrow X_{n} $ where $
\varepsilon _{n} > 0 $ and $ \varepsilon _{n} \longrightarrow 0. $
We may assume that all $ \varepsilon _{n} < \min \{ r, \varepsilon
\}. $ Then we define $ S_{n} = f_{n} (S) \subset X_{n} . $
Obviously, $ f_{n} \ : \ S \longrightarrow S_{n} $ is an $
\varepsilon _{n}-$isometry for every $ n. $ \\
Moreover, for $ x^{\prime}, \ x^{\prime\prime} \in S, $ if $
\widehat{d} ( x^{\prime}, \ x^{\prime\prime}  ) \geq \varepsilon
_{n}, $ then $ \widehat{d} ( x^{\prime}, \ x^{\prime\prime}  ) =
\widehat{d}_{n} ( f_{n} ( x^{\prime}) , \ f_{n} (
x^{\prime\prime}) ) $ (via the condition $ (SI_{2}) $ of
Def.2.22).  \\
Furthermore, let $ x \in S $ and $ y \in S_{n} $ satisfy $
\widehat{d}_{n} ( y , f_{n} ( x )) > \varepsilon _{n}. $ Firstly,
we have $ x = x_{i} $ and $ y = f_{n} ( x_{j} ) $ for some $ i, \
j \in \{ 1, \cdots , N \} $ with $ i \neq j. $ Then $ \widehat{d}
( x_{i}, \ x_{j} ) \geq r > \varepsilon _{n} . $ So by the
condition $ (SI_{2}) $ of Def.2.22 for $ f_{n} \ : \ X
\longrightarrow X_{n}, $ we have $$ \widehat{d} ( x_{i}, \ x_{j} )
= \widehat{d}_{n} ( f_{n} ( x_{i}) , \ f_{n} ( x_{j}) ) =
\widehat{d}_{n} ( y , \ f_{n} ( x) ).   $$ also we have $
\widehat{d}_{n} ( y , \ f_{n} ( x_{j}) ) =  \widehat{d}_{n} (
f_{n} ( x_{j}) , \ f_{n} ( x_{j}) ) = 0 < \varepsilon _{n}. $ So
(via taking $ x^{\prime} = x_{j} $) the condition $ ( SI_{1}) $ of
Def.2.22 holds for $ f_{n} \ : \ S \longrightarrow S_{n}. $ This
shows that $ f_{n} \ : \ S \longrightarrow S_{n} $ is a $ strong \
\varepsilon _{n}-$isometry for every $ n. $ Hence by Criterion
3.3, we obtain that $  S_{n}
\longrightarrow _{\text{GH}_{S}} S . $ \\
Lastly we come to prove that $ S_{n} $ is an $ \varepsilon -$net
in $ X_{n} $ for every $ n . $ We may as well assume that $
\varepsilon _{n} < \frac{\varepsilon }{2} . $ For every $ x \in
X_{n} $ since $ f_{n} ( X ) $ is an $ \varepsilon _{n}-$net in $
X_{n}, \ \text{dist} ( x , f_{n} ( X ) ) < \varepsilon _{n}, $ so
there exists an $ x^{\prime} \in X $ such that $ \widehat{d}_{n} (
x , f_{n} ( x^{\prime} ) ) < \varepsilon _{n}. $ Similarly, $
\widehat{d} ( x^{\prime},  x^{\prime\prime}) < \varepsilon $ for
some $ x^{\prime\prime} \in S $ because $ S $ is an $ \varepsilon
-$net in $ X . $ Then we discuss into the following two cases: \\
Case A. \ If $ \widehat{d} ( x^{\prime},  x^{\prime\prime}) \geq
\varepsilon _{n}, $ then by the condition $ (SI_{2}) $ of Def.2.22
for $ f_{n} \ : \ X \longrightarrow X_{n} , $ we have $
\widehat{d}_{n} ( f_{n} ( x^{\prime}), \  f_{n} ( x^{\prime\prime}
) ) = \widehat{d} ( x^{\prime}, x^{\prime\prime}) < \varepsilon .
$ Hence $$ \widehat{d}_{n} ( x , \  f_{n} ( x^{\prime\prime} ) )
\leq \max \{ \widehat{d}_{n} ( x , \  f_{n} ( x^{\prime} ) ) , \
\widehat{d}_{n} ( f_{n} ( x^{\prime} ), \ f_{n} ( x^{\prime\prime}
) ) \} < \varepsilon . $$ Since $ f_{n} ( x^{\prime\prime} ) \in
f_{n} ( S ) = S_{n}, $ we get $ \text{dist} ( x , S_{n} ) \leq
\widehat{d}_{n} ( x , \  f_{n} ( x^{\prime\prime} ) ) <
\varepsilon . $ \\
Case B. \ If $ \widehat{d} ( x^{\prime},  x^{\prime\prime}) <
\varepsilon _{n}, $ then since
\begin{align*}
&\mid \widehat{d} ( x^{\prime} , \  x^{\prime\prime} ) \ - \
\widehat{d}_{n} ( f_{n} ( x^{\prime} ) , \  f_{n} (
x^{\prime\prime} ) ) \mid \ \leq \text{dis} f_{n} < \varepsilon
_{n}, \ \text{we have} \\
&\widehat{d}_{n} ( f_{n} ( x^{\prime} ) , \  f_{n} (
x^{\prime\prime} )) < \varepsilon _{n} + \widehat{d} ( x^{\prime}
, \  x^{\prime\prime} ) < 2 \varepsilon _{n} < \varepsilon .
\end{align*}
$$ \text{So} \quad \widehat{d}_{n} ( x , \  f_{n} (
x^{\prime\prime} ) ) \leq \max \{ \widehat{d}_{n} ( x , \  f_{n} (
x^{\prime} ) ) , \ \widehat{d}_{n} ( f_{n} ( x^{\prime} ) , \
f_{n} ( x^{\prime\prime} ) )  \} < \varepsilon . $$ Since $  f_{n}
( x^{\prime\prime} ) \in f_{n} ( S ) = S_{n} , $ we get $
\text{dist} ( x, \ S_{n} ) \leq \widehat{d}_{n} ( x , \  f_{n} (
x^{\prime\prime} ) ) < \varepsilon . $ \\
Thus Case A and Case B together show that $ \text{dist} ( x, \
S_{n} ) < \varepsilon $ for all $ x \in X_{n}. $ So $ S_{n} $ is
an $ \varepsilon -$net in $ X_{n}. $ This completes the proof of
Theorem 3.8. \quad $ \Box $
\par  \vskip 0.2 cm

{\bf Theorem 3.9.} \ Let $ \{ X_{n} \}_{n = 1}^{\infty } $ be a
sequence of non-Archimedean metric spaces, and let $ X = \{ x_{i}
\ : \ 1 \leq i \leq N \} $ be a finite non-Archimedean metric
space of cardinality $ N. $ \ Then  \ $ X_{n} \longrightarrow
_{\text{GH}_{S}} X $ \ if and only if the
following holds: \\
For all sufficiently large $ n, \ X_{n} $ can be split into a
disjoint union of $ N $ non-empty sets $ X_{n , 1}, \ X_{n , 2}, \
\cdots , \ X_{n , N } $ so that for all $ i, \ j, $ \\
$ \text{diam} ( X_{n , i} ) \longrightarrow 0 , \quad \text{dist}
( X_{n , i}, \  X_{n , j} ) = \widehat{d} ( x_{i}, \ x_{j} ) \quad
( n \longrightarrow \infty ), $ \\
where $ \widehat{d} $ is the metric on $ X $ and $ \widehat{d}_{n}
$ is the metric on $ X_{n} $ for each $ n. $
\par  \vskip 0.2 cm

{\bf Proof.} \ Denote $ r_{0} = \min \{ \widehat{d} ( x_{i}, x_{j}
) \ : \ 1 \leq i , \ j \leq N  \ \text{and} \ i \neq j \}. $ Then
$ r_{0} > 0 . $ \\
$ \Longrightarrow : $ \ We assume that $ X_{n} \longrightarrow
_{\text{GH}_{S}} X. $ Then for any $ 0 < \varepsilon <
\frac{r_{0}}{2}, $ there exists a $ n_{0} ( \varepsilon ) \in \N $
such that $ \widehat{d}_{GH} ( X_{n}, X ) < \varepsilon $ for all
$ n \geq n_{0} ( \varepsilon ). $ By Theorem 2.23.(1), there
exists a strong $ \varepsilon -$isometry $ f_{n , \ \varepsilon }
\ : \ X_{n} \longrightarrow X $ for every $  n \geq n_{0} (
\varepsilon ). $ We claim that each $ f_{n , \ \varepsilon } $ is
surjective. If otherwise, then we may as well assume that $ x_{N}
\notin f_{n , \ \varepsilon } ( X_{n} ). $ Then $ \text{dist} (
x_{N}, \ f_{n , \ \varepsilon } ( X_{n} ) ) \geq r_{0} >
\varepsilon , $ so $ f_{n , \ \varepsilon } ( X_{n} ) $ is not an
$ \varepsilon -$net  in $ X . $ A contradiction! Therefore $ f_{n
, \ \varepsilon } $ must be surjective for every $  n \geq n_{0} (
\varepsilon ). $ Denote $$ X_{n, \ i, \ \varepsilon } = f_{n , \
\varepsilon }^{-1} ( x_{i} ) = \{ x \in  X_{n} \ : \ f_{n , \
\varepsilon } ( x ) = x_{i} \}. $$ Then $ X_{n, \ i, \ \varepsilon
} \neq \emptyset $ and $ X_{n} = \sqcup _{ 1 \leq i \leq N } X_{n,
\ i, \ \varepsilon } $ is the disjoint union of all $ X_{n, \ i, \
\varepsilon } \ ( 1 \leq i \leq N ). $ Since $ f_{n , \
\varepsilon } ( X_{n, \ i, \ \varepsilon } ) = \{ x_{i} \} $ for
each $ i , $ for any $ x^{\prime}, \ x^{\prime\prime} \in X_{n, \
i, \ \varepsilon } , $ we have $$ \widehat{d}_{n} ( x^{\prime}, \
x^{\prime\prime} ) = \ \mid  \widehat{d}_{n} ( x^{\prime}, \
x^{\prime\prime} ) \ - \ \widehat{d} ( f_{n , \ \varepsilon } (
x^{\prime} ) , \ f_{n , \ \varepsilon } ( x^{\prime\prime} )) \mid
\ \leq \text{dis} f_{n, \ \varepsilon} < \varepsilon , $$ so $
\text{diam} ( X_{n, \ i, \ \varepsilon } ) \leq \text{dis} f_{n, \
\varepsilon} < \varepsilon $ for all $ n \geq n_{0} ( \varepsilon
). $ \\
Furthermore, let $ y^{\prime} \in X_{n, \ i, \ \varepsilon }, \
y^{\prime\prime} \in X_{n, \ j, \ \varepsilon } \ ( i \neq j ), $
then  $ f_{n, \ \varepsilon} ( y^{\prime} ) = x_{i} $ \ and \ $
f_{n, \ \varepsilon} ( y^{\prime\prime} ) = x_{j} . $ Since $$
\mid  \widehat{d}_{n} ( y^{\prime}, \ y^{\prime\prime} ) \ - \
\widehat{d} ( f_{n , \ \varepsilon } ( y^{\prime} ) , \ f_{n , \
\varepsilon } ( y^{\prime\prime} )) \mid \ = \ \mid
\widehat{d}_{n} ( y^{\prime}, \ y^{\prime\prime} ) \ - \
\widehat{d} ( x_{i}, \ x_{j} ) \mid \ \leq \text{dis} f_{n, \
\varepsilon} < \varepsilon , $$ we have $ \widehat{d} ( x_{i}, \
x_{j} ) - \varepsilon < \widehat{d}_{n} ( y^{\prime}, \
y^{\prime\prime} ) < \widehat{d} ( x_{i}, \ x_{j} ) + \varepsilon
, $ so $ \widehat{d}_{n} ( y^{\prime}, \ y^{\prime\prime} ) >
\widehat{d} ( x_{i}, \ x_{j} ) - \varepsilon \geq r_{0} -
\varepsilon > 2 \varepsilon - \varepsilon = \varepsilon . $ Then
by the condition $ (SI_{2}) $ of Def.2.22, we get $
\widehat{d}_{n} ( y^{\prime}, \ y^{\prime\prime} ) = \widehat{d} (
f_{n , \ \varepsilon } ( y^{\prime} ) , \ f_{n , \ \varepsilon } (
y^{\prime\prime} )) = \widehat{d} ( x_{i}, \ x_{j} ) . $ Therefore
$$ \text{dist} ( X_{n, \ i, \ \varepsilon },
\ X_{n, \ j, \ \varepsilon } ) = \inf \{ \widehat{d}_{n} (
y^{\prime}, \ y^{\prime\prime} ) \ : \ y^{\prime} \in X_{n, \ i, \
\varepsilon }, \ y^{\prime\prime} \in X_{n, \ j, \ \varepsilon }
 \} = \widehat{d} ( x_{i}, \ x_{j} ). $$ To sum up, we have obtain
the following \\
Conclusion. \ For any $ 0 < \varepsilon  < r_{0} / 2, $ there
exists a $ n_{0}( \varepsilon ) \in \N $ such that for all $ n
\geq n_{0}( \varepsilon ), $ we have \ $ X_{n} = \sqcup _{ 1 \leq
i \leq N } X_{n, \ i, \ \varepsilon } $ ( the disjoint union )
with non-empty sets  $ X_{n, \ i, \ \varepsilon } \ ( 1 \leq i
\leq N ) $ satisfying $ \text{diam} ( X_{n, \ i, \ \varepsilon } )
< \varepsilon  $ \ and \ $ \text{dist} ( X_{n, \ i, \ \varepsilon
} , \ X_{n, \ j, \ \varepsilon } ) = \widehat{d} ( x_{i}, \ x_{j}
) . $ \\
Now we take a sequence $ \{ \varepsilon _{k} \}_{k = 1}^{\infty }
$ of positive real numbers such that $ \varepsilon _{k} < r_{0} /
2 $ for all $ k $ and $ \varepsilon _{k} \longrightarrow 0 $ as $
k \longrightarrow \infty . $ Then for each $ \varepsilon _{k} , $
we have a $ n_{0}( \varepsilon _{k}) \in \N $ such that the above
Conclusion holds for all $ n \geq n_{0}( \varepsilon _{k}). $ \\
Case A. \ If $ \{ n_{0}( \varepsilon _{k}) \}_{k = 1}^{\infty } $
is bounded, then there exists a $ m_{0} \in \N $ such that $
n_{0}( \varepsilon _{k}) < m_{0} $ for all $ k . $ So for each $
\varepsilon _{k}, $ the above Conclusion holds for all $ n \geq
m_{0}. $ Now for each $ n \geq m_{0} $ and $ 1 \leq i \leq N, $ we
take $ X_{n, \ i } = X_{n, \ i, \ \varepsilon _{n}}. $ Then we
obtain a split for all $ n \geq m_{0} $ satisfying \\
$ X_{n} = \sqcup _{ 1 \leq i \leq N } X_{n, \ i, \ \varepsilon
_{n} } = \sqcup _{ 1 \leq i \leq N } X_{n, \ i }, $ \\
$ \text{diam} ( X_{n, \ i } ) = \text{diam} ( X_{n, \ i, \
\varepsilon _{n} } ) < \varepsilon _{n} \longrightarrow 0 $ as $ n
\longrightarrow \infty , $ \  and  \\
$ \text{dist} ( X_{n, \ i } , \ X_{n, \ j } ) = \text{dist} (
X_{n, \ i, \ \varepsilon _{n} } , \ X_{n, \ j, \ \varepsilon _{n}
} ) = \widehat{d} ( x_{i}, \ x_{j} ) . $ \\
Case B. \ If $ \{ n_{0}( \varepsilon _{k}) \}_{k = 1}^{\infty } $
is unbounded, then without loss of generality, we may assume \ $
n_{0}( \varepsilon _{1}) < n_{0}( \varepsilon _{2}) < \cdots $ \
and $ n_{0}( \varepsilon _{k}) \longrightarrow \infty $ as $ k
\longrightarrow \infty . $ \ Now for each $ 1 \leq i \leq N $ and
$ n \geq n_{0}( \varepsilon _{1}), $ we take \\
$ X_{n, \ i } = X_{n, \ i, \ \varepsilon _{k} } $ \ if \ $ n_{0}(
\varepsilon _{k}) \leq n < n_{0}( \varepsilon _{k + 1}). $ \\
Then obviously we have \\
$ X_{n} = \sqcup _{ 1 \leq i \leq N }
X_{n, \ i, \ \varepsilon _{k} } = \sqcup _{ 1 \leq i \leq N } X_{n, \ i }, $ \\
$ \text{diam} ( X_{n, \ i } ) = \text{diam} ( X_{n, \ i, \
\varepsilon _{k} } ) < \varepsilon _{k} \longrightarrow 0 \quad (
n
\longrightarrow \infty ), $ \ and  \\
$ \text{dist} ( X_{n, \ i } , \ X_{n, \ j } ) = \text{dist} (
X_{n, \ i, \ \varepsilon _{k} } , \ X_{n, \ j, \ \varepsilon _{k}
} ) = \widehat{d} ( x_{i}, \ x_{j} ) . $ \\
This proves the necessity. \\
$ \Longleftarrow : $ \ Now we assume such split exists, i.e., for
all sufficiently large $ n, \ X_{n} = \sqcup _{ 1 \leq i \leq N }
X_{n, \ i } $ with the given properties. Then we define a map $
f_{n}: \ X_{n} \longrightarrow X $ by $ f_{n} (x^{\prime}) = x_{i}
$ if $ x^{\prime} \in X_{n, \ i } \quad ( i \in \{ 1, \cdots , N
\}). $ Then $ f_{n} ( X_{n, \ i } ) = \{  x_{i} \}. $ Obviously $
f_{n} $ is surjective. \\
For any $ 0 < \varepsilon < r_{0} / 2, $ there exists a $ n_{0}
\in \N $ such that $ \text{diam} ( X_{n, \ i } ) < \varepsilon $
for every $ i $ and all $ n > n_{0}. $ Then by the fact that $
\text{dist} ( X_{n, \ i }, \ X_{n, \ j } ) = \widehat{d} ( x_{i},
 x_{j} ) \geq r_{0} > \varepsilon \ ( i \neq j ) $ and the strong
 triangle inequality, we can easily obtain that
 $ \widehat{d}_{n} ( x^{\prime} , \ x^{\prime\prime} ) = \widehat{d} ( x_{i},
 x_{j} ) $ for any $ x^{\prime} \in X_{n, \ i } $ and
 $ x^{\prime\prime} \in X_{n, \ j } \ ( i \neq j ) . $ From this we have
 \begin{align*}
 \text{dis} f_{n} &= \sup \{ \mid \widehat{d}_{n} ( x^{\prime} ,
 \ x^{\prime\prime} )  \ - \ \widehat{d}( f_{n} (x^{\prime}) , \
f_{n} (x^{\prime\prime})) \mid \ : \ x^{\prime} ,
 \ x^{\prime\prime} \in X_{n} \} \\
&= \sup \{ \mid \widehat{d}_{n} ( x^{\prime} ,
 \ x^{\prime\prime} )  \mid \ : \ x^{\prime} ,
 \ x^{\prime\prime} \in X_{n , \ i} \quad ( 1 \leq i \leq N ) \} \\
 &= \max \{ \text{diam} ( X_{n , \ i} ) \ : \ 1 \leq i \leq N  \} <
 \varepsilon .
\end{align*}
Because $ f_{n} $ is surjective, so $ f_{n} $ is an $ \varepsilon
-$isometry for every $ n > n_{0}. $ \\
Now let $ x \in X_{n} $ and $ y \in X, $ then $ x \in X_{n , i } $
and $ y = f_{n} (x^{\prime}) $ with $ x^{\prime} \in X_{n , j } $
for some $ i, \ j \in \{ 1, \cdots , N \}. $ If $ \widehat{d} ( y,
f_{n} (x)) \geq \varepsilon , $ i.e., $ \widehat{d} ( f_{n}
(x^{\prime}), f_{n} (x)) \geq \varepsilon , $ then $ \widehat{d} (
x_{i} , x_{j} ) = \widehat{d} ( f_{n} (x^{\prime}), f_{n} (x))
\geq \varepsilon , $ so $ i \neq j. $ Hence as above discussed, $$
\widehat{d}_{n} ( x , \ x^{\prime} ) = \text{dist} ( X_{n , i } ,
\ X_{n , j } ) = \widehat{d} ( x_{i} , x_{j} ) = \widehat{d} ( y ,
f_{n} (x) ). $$ So the condition $ (SI_{1}) $ of Def.2.22 holds
for $ f_{n} . $ \\
Next let $ x^{\prime} , \ x^{\prime\prime} \in X_{n} , $ then $
x^{\prime} \in X_{n , i } $ and $ x^{\prime\prime} \in X_{n , j }
$ for some $ i, \ j \in \{ 1, \cdots , N \}. $ If $
\widehat{d}_{n} ( x^{\prime} , \ x^{\prime\prime} ) \geq
\varepsilon , $ then $ i \neq j. $ So $ \widehat{d}_{n} (
x^{\prime} , \ x^{\prime\prime} ) = \widehat{d}( x_{i} , \ x_{j} )
= \widehat{d}( f_{n} ( x^{\prime} ), \ f_{n} ( x^{\prime\prime})).
 $ So the condition $ ( SI_{2}) $ of Def. 2.22 holds for $ f_{n} .
 $ Therefore, for all $ n > n_{0} ,  f_{n} $ is a strong $ \varepsilon
 -$isometry from $ X_{n} $ to $ X . $ Hence by Theorem 2.23.(2),
 we obtain that $ \widehat{d}_{GH} ( X_{n} , X ) \leq \varepsilon
 $ for all $ n > n_{0}. $ This shows that $ X_{n} \longrightarrow
_{\text{GH}_{S}} X , $ and the proof of Theorem 3.9 is completed.
\quad $ \Box $
\par  \vskip 0.2 cm

{\bf Note added for the proof of Theorem 3.9.} \\
If we assume that $ X_{n} $ and $ X $ are compact, then we can
prove the sufficiency of Theorem 3.9 as
follows: \\
$ \Longleftarrow : $ \ Take $ X_{n}^{\prime } = \{ x_{1}^{\prime},
\ \cdots , \ x_{N}^{\prime} \} \subset X_{n} $ by choosing one
element $ x_{i}^{\prime} \in X_{n, \ i } $ for each $ i \in \{1, \
\cdots , \ N \}. $ Then it is easy to see that $ X_{n}^{\prime } $
is an $ \varepsilon -$net in $ X_{n} $ and $ \widehat{d}_{n}
(x_{i}^{\prime}, \ x_{j}^{\prime} ) = \widehat{d} ( x_{i}, \
x_{j}) $ for all $ i, j. $ So $ X_{n} $ and $ X $ are $ strong \
\varepsilon -approximations $ of each other. Therefore by Theorem
3.5.(1), we have $ \widehat{d}_{GH} ( X_{n} , X ) \leq \varepsilon
. $ \quad $ \Box $
\par  \vskip 0.2 cm

{\bf Definition 3.10.} \ Let $ ( X, d_{X} ) $ be a metric space
and $ S $ be a non-empty subset of $ X. $ We denote $$ W_{X}(S) =
\{ d_{X} ( x, \ y ) : \ x, \ y \in S \ \text{and} \ x \neq y \}
\subset \R_{> 0 } \cup \{ \infty \}, $$ and write $ \omega _{0} =
\inf \{ \alpha : \ \alpha \in W_{X}(S) \}, \quad \omega _{1} =
\sup \{ \alpha : \ \alpha \in W_{X}(S) \}. $ We call $ W_{X}(S) $
the metric weight set associated to $ S $ in $ X. $
\par  \vskip 0.2 cm

{\bf Definition 3.11.} \ We say that a class $
\widehat{\mathcal{X}}_{sut} $ of compact non-Archimedean metric
spaces is $ strongly \ uniformly \ totally \ bounded $ \ if the
following holds: \\
For every $ \varepsilon > 0, $ there exists a positive integer
number $ N = N ( \varepsilon ) $ and a finite set $ R(\varepsilon)
\subset \R_{> 0} $ of positive real numbers such that every $ X
\in \widehat{\mathcal{X}}_{sut} $ contains an $ \varepsilon -$net
$ S_{(X)} $ consisting of no more than $ N $ points and $
W_{X}(S_{(X)}) \subset R(\varepsilon) . $
\par  \vskip 0.2 cm

Now we establish the following compactness Theorem for the class $
\widehat{\mathcal{X}}_{sut} $ under $ \widehat{d}_{GH}. $ For the
corresponding theorem about $ uniformly \ totally \ bounded $
class of compact metric spaces under $ d_{GH}, $ see [BBI,
Thm.7.4.15].
\par  \vskip 0.2 cm

{\bf Theorem 3.12} ( Compactness Theorem ). \\
Any $ strongly \ uniformly \ totally \ bounded $ class $
\widehat{\mathcal{X}}_{sut} $ of compact non-Archimedean metric
spaces is pre-compact in the strong Gromov-Hausdorff topology.
That is, any sequence of elements of $ \widehat{\mathcal{X}}_{sut}
$ contains a Cauchy subsequence under the metric $
\widehat{d}_{GH}. $
\par  \vskip 0.2 cm

{\bf Proof.} \ Let $ \{ X_{n} \}_{n = 1}^{\infty } $ be a sequence
in $ \widehat{\mathcal{X}}_{sut} , $ we denote by $
\widehat{d}_{n} $ the metric of $ X_{n} . $ In every space $
X_{n}, $ there exists an $ 1-$net $ T_{n}^{(1)} = S_{n}^{(1)} $
consisting of no more than $ N_{1} = N_{1} (1) $ points, and $
W_{X_{n}}(S_{n}^{(1)}) \subset R(1). $ Let $ n_{1} = \sharp
S_{n}^{(1)} \leq  N_{1}, $ and write
\begin{align*}
S_{n}^{(1)} = & \{ x_{n, \ 1}, \ \cdots , \ x_{n, \ n_{1} }, \
x_{n, \ n_{1} + 1 }, \ \cdots , \ x_{n, \ N_{1} } \} \quad
\text{with} \\
& x_{n, \ n_{1} } = x_{n, \ n_{1} + 1 } = \cdots  = x_{n, \ N_{1}
}.
\end{align*}
Likewise, $ X_{n} $ contains a $ 1/2-$net $ S_{n}^{(2)} $ with $
\sharp S_{n}^{(2)} = n_{2} \leq  N (1/2) = N_{2} $ and $
W_{X_{n}}(S_{n}^{(2)}) \subset R(1/2). $ Denote $ T_{n}^{(2)} =
T_{n}^{(1)} \cup S_{n}^{(2)} $ and write
\begin{align*}
&T_{n}^{(2)} \\
&= \{ x_{n, \ 1}, \ \cdots , \ x_{n, \ N_{1} }, \ x_{n, \ N_{1} +
1 }, \ \cdots , \ x_{n, \ N_{1} + n_{2} }, \ x_{n, \ N_{1} + n_{2}
+ 1 }, \ \cdots , \ x_{n, \ N_{1} + N_{2} } \} \\
&\text{with} \quad S_{n}^{(2)} = \{ x_{n, \ N_{1} + 1 }, \ \cdots
, \ x_{n, \ N_{1} + n_{2} } \} \\
&\text{and} \quad  x_{n, \ N_{1} + n_{2} } = x_{n, \ N_{1} + n_{2}
+ 1 } = \cdots = x_{n, \ N_{1} + N_{2} }.
\end{align*}
Obviously, $ T_{n}^{(2)} $ is also a $ 1/2-$net in $ X_{n}. $
Follows this way, we obtain that
\begin{align*}
&T_{n}^{(k)} \\
&= \{ x_{n, \ 1}, \ \cdots , \ x_{n, \ L_{k - 1} }, \ x_{n, \ L_{k
- 1} + 1 }, \ \cdots , \ x_{n, \ L_{k - 1} + n_{k} }, \ x_{n, \
L_{k - 1} + n_{k} + 1 }, \ \cdots , \ x_{n, \ L_{k - 1} + N_{k} } \}, \\
&\text{where} \quad  x_{n, \ L_{k - 1} + n_{k} } = x_{n, \ L_{k -
1} + n_{k} + 1 } = \cdots = x_{n, \ L_{k - 1} + N_{k} }, \\
&\text{and} \quad S_{n}^{(k)} = \{ x_{n, \ L_{k - 1} + 1 }, \
\cdots , \ x_{n, \ L_{k - 1} + n_{k} } \}
\end{align*}
is a $ 1 / k -$net in $ X_{n} $ with $ n_{k} = \sharp S_{n}^{(k)}
 \leq  N (1 / k ) = N_{k}, $ and $ W_{X_{n}}(S_{n}^{(k)}) \subset
R(1 / k ), $ and $ L_{k} = N_{1} + N_{2} + \cdots + N_{k}. $ \\
Obviously, $ T_{n}^{(k)} $ is a $ 1 / k -$net in $ X_{n} $ and $
T_{n}^{(k)} = T_{n}^{(k - 1 )} \cup S_{n}^{(k)} = S_{n}^{(1)} \cup
S_{n}^{(2)} \cup \cdots \cup S_{n}^{(k)}. $ Let $ k
\longrightarrow \infty , $ we obtain a countable dense subset $
T_{n} =  \{ x_{n, \ i} \}_{i = 1}^{\infty } \subset X_{n} $ such
that for every $ k, $ the first $ L_{k} $ points of $ T_{n} $ form
a $ 1 / k -$net in $ X_{n} $ ( some points in $ T_{n} $ may
coincide ). The density is easy to see, since for every $ x \in
X_{n} $ and $ \varepsilon > 0 , $ there exists $ k \in \N $ such
that $ \varepsilon > 1 / k $ and then $ \text{dist} ( x, \
S_{n}^{(k)} ) < 1 / k < \varepsilon , $ hence $ \text{dist} ( x, \
T_{n} ) <  \varepsilon $ because $ S_{n}^{(k)} \subset T_{n}. $ \\
Denote $ D = \max \{ 1, \ \max \{ a : \ a \in R(1) \} \}. $ Then
for any $ x^{\prime}, \ x^{\prime\prime} \in X_{n}, $ we have $
\text{dist}( x^{\prime}, \ S_{n}^{(1)} ) < 1 $ and $ \text{dist}(
x^{\prime\prime}, \ S_{n}^{(1)} ) < 1 , $ so there exist $ y_{1},
\ y_{2} \in S_{n}^{(1)} $ such that $ \widehat{d}_{n} (
x^{\prime}, \ y_{1} ) < 1 $ and $ \widehat{d}_{n} (
x^{\prime\prime}, \ y_{2} ) < 1. $ Then $$ \widehat{d}_{n} (
x^{\prime}, \ x^{\prime\prime} ) \leq \max \{ \widehat{d}_{n} (
x^{\prime}, \ y_{1} ), \ \widehat{d}_{n} ( y_{1}, \ y_{2} ), \
\widehat{d}_{n} ( y_{2}, \ x^{\prime\prime} ) \} \leq D $$ because
$ W_{X_{n}}(S_{n}^{(1)}) \subset R(1). $ This is independent of $
n. $ Therefore, for each pair $ ( i, \ j ), \ \widehat{d}_{n} (
x_{n , \ i}, \ x_{n , \ j} ) \leq D, $ i.e., $ \ \widehat{d}_{n} (
x_{n , \ i}, \ x_{n , \ j} ) \in [0, \ D], $ the compact interval.
Hence the sequence $ \{ \widehat{d}_{n} ( x_{n , \ i}, \ x_{n , \
j} ) \}_{n = 1}^{\infty } $ contains a converging subsequence.
Then using the Cantor diagonal procedure, we can extract a
subsequence of $ \{ X_{n} \}_{n = 1}^{\infty } $ in which $ \{
\widehat{d}_{n} ( x_{n , \ i}, \ x_{n , \ j} ) \}_{n = 1}^{\infty
} $ converges for all $ i, \ j. $ Without loss of generality, we
may assume that they converge without passing to a subsequence.
Now we come to construct a complete non-Archimedean metric space $
\overline{X} $ as follows: \\
Pick an abstract countable set $ X =  \{ x_{i} \}_{i = 1}^{\infty
} $ and define a function $ \widehat{d} $ on $ X \times X $ by  $$
\widehat{d}( x_{i}, \ x_{j} ) = \lim _{ n \longrightarrow \infty }
\widehat{d}_{n} ( x_{n , \ i}, \ x_{n , \ j} ) \quad ( \forall \
i, \ j \in \N ). $$ For any $ x_{k}, $ since
\begin{align*}
&\widehat{d}_{n} ( x_{n , \ i}, \ x_{n , \ j} ) \leq \max \{
\widehat{d}_{n} ( x_{n , \ i}, \ x_{n ,
\ k} ) , \ \widehat{d}_{n} ( x_{n , \ k}, \ x_{n , \ j} ) \} \\
&=\frac{1}{2} ( \widehat{d}_{n} ( x_{n , \ i} , \ x_{n , \ k} )  +
\widehat{d}_{n} ( x_{n , \ k}, \ x_{n , \ j} )) + \frac{1}{2} \mid
\widehat{d}_{n} ( x_{n , \ i} , \ x_{n , \ k} )  - \widehat{d}_{n}
( x_{n , \ k}, \ x_{n , \ j} ) \mid  ,
\end{align*}
we have
\begin{align*}
\widehat{d}( x_{i}, \ x_{j} ) & \leq \frac{1}{2} ( \lim _{ n
\longrightarrow \infty } \widehat{d}_{n} ( x_{n , \ i} , \ x_{n ,
\ k} )  + \lim _{ n \longrightarrow \infty } \widehat{d}_{n} (
x_{n , \ k}, \ x_{n , \ j} )) \\
& + \frac{1}{2} \mid \lim _{ n \longrightarrow \infty }
\widehat{d}_{n} ( x_{n , \ i} , \ x_{n , \ k} )  - \lim _{ n
\longrightarrow \infty } \widehat{d}_{n} ( x_{n , \ k}, \ x_{n , \
j} ) \mid \\
& = \frac{1}{2} ( \widehat{d} ( x_{i} , \ x_{k} ) + \widehat{d} (
x_{k}, \ x_{j} )) + \frac{1}{2} \mid \widehat{d}( x_{i} , \ x_{k}
)  - \widehat{d}( x_{k}, \ x_{j} ) \mid \\
& = \max \{ \widehat{d}( x_{i}, \ x_{k} ) , \ \widehat{d}( x_{k},
\ x_{j} ) \}.
\end{align*}
So $ \widehat{d} $ satisfies the strong triangle inequality. The
symmetry is obvious, so $ \widehat{d} $ is a non-Archimedean
semi-metric on $ X, $ and then the quotient space $ X /
\widehat{d} $ is a non-Archimedean metric space. We will denote by
$ \overline{x}_{i} $ the point of $ X / \widehat{d} $ obtained
from $ x_{i}. \ X / \widehat{d} $ may not be complete, so let $
\overline{X} $ be the completion of $ X / \widehat{d}. $ \\
Let $ \varepsilon > 0 $ and take a $ k \in \N $ such that $ 1 / k
< \varepsilon . $ We consider the set $ \overline{T}^{(k)} = \{
\overline{x}_{i} : \ 1 \leq i \leq L_{k} \} $ with $ L_{k} $ as
above. For any $ \overline{x} \in X / \widehat{d} , \ \overline{x}
= \overline{x}_{j} $ for some $ j \in \N. $ Then for $ x_{n , \ j}
\in T_{n} \subset X_{n}, $ \ since $ T_{n}^{(k)} $ is a $ 1 / k
-$net in $ X_{n}, $ there exist some $ l \leq L_{k} $ such that $
\widehat{d}_{n} ( x_{n , \ j}, \ x_{n , \ l} ) < 1 / k . $ Since $
L_{k} $ is independent of $ n, $ there exists a $ l \leq L_{k} $
such that $ \widehat{d}_{n} ( x_{n , \ j}, \ x_{n , \ l} ) < 1 / k
$ for infinitely many indices $ n. $ Then $ \widehat{d}( x_{j}, \
x_{l} ) = \lim _{ n \longrightarrow \infty } \widehat{d}_{n} (
x_{n , \ j}, \ x_{n , \ l} ) \leq 1 / k $ ( the limit of a
converging sequence is equal to the  limit of its any subsequence
). So $ \widehat{d} ( \overline{x}, \overline{x}_{l}) =
\widehat{d} ( \overline{x}_{j}, \overline{x}_{l}) = \widehat{d} (
x_{j}, x_{l}) < \varepsilon . $ This shows that $
\overline{T}^{(k)} $ is an $ \varepsilon -$net in $ X /
\widehat{d}. $ Since $ X / \widehat{d} $ is dense in $
\overline{X}, $ for any $ \overline{x} \in \overline{X}, $ there
exists an $ \overline{x}^{\prime} \in X / \widehat{d} $ such that
$ \widehat{d}( \overline{x}, \ \overline{x}^{\prime} ) <
\varepsilon . $ For such $ \overline{x}^{\prime} , \ \text{dist} (
\overline{x}^{\prime}, \ \overline{T}^{(k)} ) < \varepsilon , $ so
there exists an $ \overline{x}_{i} \in \overline{T}^{(k)} $ such
that $ \widehat{d}( \overline{x}^{\prime}, \ \overline{x}_{i} ) <
\varepsilon . $ Hence $ \widehat{d}( \overline{x}, \
\overline{x}_{i} ) \leq \max \{ \widehat{d}( \overline{x}, \
\overline{x}^{\prime} ) , \ \widehat{d}( \overline{x}^{\prime}, \
\overline{x}_{i} ) \} < \varepsilon , $ which shows that $
\overline{T}^{(k)} $ is an $ \varepsilon -$net in $ \overline{X}.
$ Thus $ \overline{X} $ is totally bounded because $ \varepsilon
> 0 $ is arbitrary, hence $ \overline{X} $ is compact since it is
complete. Therefore, $ (\overline{X}, \ \widehat{d} ) $ is a
compact non-Archimedean metric space, i.e., its isometry class
belongs to $ \widehat{\Gamma }_{c}. $ \\
Furthermore, for $ \varepsilon > 0 $ and $ k \in \N $ with $ 1 / k
< \varepsilon $ as above, we consider the sets $ S_{n}^{(k)} = \{
x_{n, \ L_{k - 1} + 1 }, \ \cdots , \ x_{n, \ L_{k - 1} + n_{k}}
\} \subset T_{n}^{k} $ for every $ n $ as before. Since $ N_{k} =
N( 1 / k ) $ is independent of $ n $ and $ n_{k} = \sharp
S_{n}^{(k)} \leq N_{k}, $ there exists a positive integer $
N^{\prime} \leq N_{k} $ such that $ n_{k} = N^{\prime} $ for
infinitely many $ n. $ So we may as well assume that $ n_{k}
= N^{\prime} $ for all $ n . $ \\
Without loss of generality, we may assume that $ \sharp R ( 1 / k
) > 1 , $ and then define \\
$ r_{k} = \min \{ \mid a - b  \mid \
: \ a, \ b \in R ( 1 / k ) \ \text{and} \ a \neq b \} . $
Obviously, $ r_{k} > 0 . $ \\
Define \ $ \overline{S}^{(k)} = \{ \overline{x}_{i} : \ L_{k - 1}
+ 1 \leq i \leq L_{k - 1} + N^{\prime} \} \subset \overline{X}, $
where $ L_{k - 1} = N_{1} + N_{2} + \cdots + N_{k - 1} $ is
independent of $ n $ as above. For any $ \overline{x} \in X /
\widehat{d} , \ \overline{x} = \overline{x}_{j} $ for some $ j \in
\N . $ Since $ S_{n}^{(k)} $ is a $ 1 / k -$net in $ X_{n} $ for
every $ n , $ and $ x_{n, \ j} \in T_{n} \subset X_{n}, $ there
exist some $ l_{n} : \ L_{k - 1} + 1 \leq  l_{n} \leq L_{k - 1} +
N^{\prime} $ such that $ \widehat{d}_{n} ( x_{n, \ j}, \ x_{n, \
l_{n} }) < 1 / k , $ and then there exist at least one $ l $ such
that $ L_{k - 1} + 1 \leq  l \leq L_{k - 1} + N^{\prime} $ and $
\widehat{d}_{n} ( x_{n, \ j}, \ x_{n, \ l }) < 1 / k $ for
infinitely many indices $ n. $ Thus
\begin{align*}
&\widehat{d}( x_{j}, \ x_{l }) = \lim _{n \longrightarrow \infty }
\widehat{d}_{n} ( x_{n, \ j}, \ x_{n, \ l }) \leq \frac{1}{k} <
\varepsilon , \quad \text{so} \\
&\widehat{d}( \overline{x}, \ \overline{x}_{l} ) = \widehat{d}(
\overline{x}_{j}, \ \overline{x}_{l} ) = \widehat{d}( x_{j}, \
x_{l }) < \varepsilon , \quad \text{ and then} \\
&\text{dist} ( \overline{x}, \ \overline{S}^{(k)}) \leq
\widehat{d}( \overline{x}, \ \overline{x}_{l} ) < \varepsilon .
\end{align*}
Hence $ \overline{S}^{(k)} $ is an $ \varepsilon -$net in $ X /
\widehat{d}. $ Because $ X / \widehat{d} $ is dense in $
\overline{X}, $ like the above proof for $ \overline{T}^{(k)}, $
one can show that $ \overline{S}^{(k)} $ is also an $ \varepsilon
-$net in $ \overline{X}. $ \\
Now for any $ i, \ j : \ L_{k - 1} + 1 \leq i, \ j \leq L_{k - 1}
+ N^{\prime} $ and $ i \neq j , $ we have as before $ \widehat{d}(
\overline{x}_{i}, \ \overline{x}_{j}) = \widehat{d}( x_{i}, \ x_{j
}) = \lim _{n \longrightarrow \infty } \widehat{d}_{n} ( x_{n, \
i}, \ x_{n, \ j }) $ with $ x_{n, \ i}, \ x_{n, \ j } \in
S_{n}^{(k)} . $ Then for any $ 0 < \delta < r_{k}, $ there exists
$ n_{0} \in \N $ such that for all $ n_{1}, \ n_{2} > n_{0} $  we
have $$ \mid \widehat{d}_{n_{1}} ( x_{n_{1}, \ i}, \ x_{n_{1}, \ j
}) - \widehat{d}_{n_{2}} ( x_{n_{2}, \ i}, \ x_{n_{2}, \ j }) \mid
\ < \delta . $$ Since $ \widehat{d}_{n_{1}} ( x_{n_{1}, \ i}, \
x_{n_{1}, \ j }) \in W_{X_{n_{1}}} ( S_{n_{1}}^{(k)}) \subset R(1
/ k) $ and $ \widehat{d}_{n_{2}} ( x_{n_{2}, \ i}, \ x_{n_{2}, \ j
}) \in W_{X_{n_{2}}} ( S_{n_{2}}^{(k)}) \subset R(1 / k), $ we
obtain that $$ \widehat{d}_{n_{1}} ( x_{n_{1}, \ i}, \ x_{n_{1}, \
j }) = \widehat{d}_{n_{2}} ( x_{n_{2}, \ i}, \ x_{n_{2}, \ j }) \
\text{or} \ \mid \widehat{d}_{n_{1}} ( x_{n_{1}, \ i}, \ x_{n_{1},
\ j }) - \widehat{d}_{n_{2}} ( x_{n_{2}, \ i}, \ x_{n_{2}, \ j })
\mid \
\geq r_{k}. $$ But $ \delta < r_{k}, $ so we must have \\
$ \widehat{d}_{n_{1}} ( x_{n_{1}, \ i}, \ x_{n_{1}, \ j }) =
\widehat{d}_{n_{2}} ( x_{n_{2}, \ i}, \ x_{n_{2}, \ j }) $ for all
$ n_{1}, \ n_{2} > n_{0}. $ This shows that \\
$ \widehat{d}_{n} ( x_{n, \ i}, \ x_{n, \ j }) = \widehat{d} (
\overline{x}_{i}, \ \overline{x}_{j} ) $ for all sufficiently
large $ n . $ Since the number of such pairs $ ( i , \ j ) $ is
finite, it follows that, for all sufficiently large $ n , \
\widehat{d}_{n} ( x_{n, \ i}, \ x_{n, \ j }) = \widehat{d} (
\overline{x}_{i}, \ \overline{x}_{j} ) $ for all $ i, \ j $
satisfying $ L_{k - 1 } + 1 \leq i, \ j \leq L_{k - 1 } +
N^{\prime} . $ Therefore, by Theorem 3.9, we get $ S_{n}^{(k)}
\longrightarrow _{\text{GH}_{S}} \overline{S}^{(k)}. $ Since $
S_{n}^{(k)} $ is an $ \varepsilon -$net in $ X_{n} $ for every $
n, $ and $ \overline{S}^{(k)} $ is an $ \varepsilon -$net in $
\overline{X}, $ by Theorem 3.8, we obtain that $ X_{n}
\longrightarrow _{\text{GH}_{S}} \overline{X}. $ This completes
the proof of Theorem 3.12. \quad $ \Box $

\par     \vskip  1 cm

\hspace{-0.6cm}{\bf 4. \ Computing the Non-Archimedean
Gromov-Hausdorff distance }

\par \vskip 0.8 cm

{\bf Definition 4.1.} \ Let $ ( X, d_{X}) $ be a metric space and
$ S $ be a non-empty subset of $ X, \ W_{X} ( S ) $ is the metric
weight set associated to $ S $ in $ X $ as before.  Let $
\varepsilon
> 0, $ we define $$ W_{X} ( S )_{ \geq \varepsilon} = \{
 d_{X}( x, \ y ): \ x, \ y \in S \ \text{and} \ d_{X} ( x, \ y )
 \geq \varepsilon \} = W_{X} ( S ) \cap [\varepsilon , \ \infty ]. $$
\par  \vskip 0.2 cm

{\bf Theorem 4.2.} \ Let $ ( X, \ \widehat{d}_{X}) $ and $ ( Y, \
\widehat{d}_{Y}) $ be two non-Archimedean metric spaces and denote
$ D = \max \{ \text{diam}(X), \ \text{diam} (Y) \}. $ Then
\par  \vskip 0.15 cm
(1) \ $ \widehat{d}_{GH} ( X, \ Y ) \geq \inf \{ \varepsilon > 0 :
\ W_{X} ( X )_{ \geq \varepsilon } = W_{Y} ( Y )_{ \geq
\varepsilon } \} . $
\par  \vskip 0.15 cm
(2) \ If there does not exist $ \varepsilon > 0 $ such that $ \
W_{X} ( X )_{ \geq \varepsilon } =
 W_{Y} ( Y )_{ \geq \varepsilon } , $
then $ \widehat{d}_{GH} ( X, \ Y ) = \infty . $
\par  \vskip 0.15 cm
(3A) \ If $ D < + \infty , $ then $  \widehat{d}_{GH} ( X, \ Y )
\leq D. $
\par  \vskip 0.15 cm
(3B) \ If $ D < + \infty $ and $ \text{diam} (X) \neq \text{diam}
(Y) , $ then $ \widehat{d}_{GH} ( X, \ Y ) = D . $
\par  \vskip 0.2 cm

{\bf Proof.} \ We denote $ d = \inf \{ \varepsilon > 0 : \ W_{X} (
X )_{ \geq \varepsilon} = W_{Y} ( Y )_{ \geq \varepsilon} \} . $
\par  \vskip 0.15 cm
(1) \ If $ \widehat{d}_{GH} ( X, \ Y ) = + \infty , $ then we are
done. So we assume that $ \widehat{d}_{GH} ( X, \ Y ) < + \infty .
$ For any $ \varepsilon > \widehat{d}_{GH} ( X, \ Y ) , $ we need
to prove that $ W_{X} ( X )_{ \geq \varepsilon } = W_{Y} ( Y )_{
\geq \varepsilon } . $ \\
If both $ W_{X} ( X )_{ \geq \varepsilon} = \emptyset  $ and $
W_{Y} ( Y )_{ \geq \varepsilon} = \emptyset , $ then we are done.
So we may as well assume that $ W_{X} ( X )_{ \geq \varepsilon }
\neq \emptyset . $ By Theorem 2.23.(1), there exists a $ strong \
\varepsilon -$isometry $ f : \ X \longrightarrow Y. $ Let $ r \in
W_{X} ( X )_{ \geq \varepsilon }, $ then $ r = \widehat{d}_{X} (
x_{1},  x_{2} ) \geq \varepsilon $ for some $ x_{1}, \ x_{2} \in
X. $ So by the condition $ (SI_{2}) $ of Def.2.22, we have $ r =
\widehat{d}_{Y} ( f(x_{1}),  f(x_{2}) ) \in W_{Y} ( Y )_{ \geq
\varepsilon }. $ This shows that $ W_{X} ( X )_{ \geq \varepsilon
} \subset W_{Y} ( Y )_{ \geq \varepsilon }, $ in particular $
W_{Y} ( Y )_{ \geq \varepsilon } \neq \emptyset . $ \\
From $ \widehat{d}_{GH} ( Y, \ X ) = \widehat{d}_{GH} ( X, \ Y ) <
\varepsilon , $ there also exists by Theorem 2.23.(1) a $ strong \
\varepsilon -$isometry $ g : \ Y \longrightarrow X. $ Let $
r^{\prime} \in W_{Y} ( Y )_{ \geq \varepsilon }, $ then $
r^{\prime} = \widehat{d}_{Y}( y_{1}, \ y_{2} ) \geq \varepsilon $
for some $ y_{1}, \ y_{2} \in Y. $ So by the condition $ (SI_{2})
$ of Def.2.22, we have $ r^{\prime} = \widehat{d}_{X} ( g(y_{1}),
g(y_{2}) ) \in W_{X} ( X )_{ \geq \varepsilon }. $ This shows that
$ W_{Y} ( Y )_{ \geq \varepsilon } \subset W_{X} ( X )_{ \geq
\varepsilon }. $ Hence $ W_{X} ( X )_{ \geq \varepsilon } = W_{Y}
( Y )_{ \geq \varepsilon }. $ Therefore, $ d \leq \varepsilon . $
Since $ \varepsilon > \widehat{d}_{GH} ( X, \ Y ) $ is arbitrary,
we get $ d \leq \widehat{d}_{GH} ( X, \ Y ) . $ This proves (1).
\par  \vskip 0.15 cm
(2) \ Easily follows from (1).
\par  \vskip 0.15 cm
(3A) \ We define a function $ \widehat{d} $ on $ Z \times Z $ with
$ Z = X \sqcup Y $ (disjoint
union ) as follows: \\
$ \widehat{d} \mid _{X \times X} = \widehat{d}_{X} ; \ \widehat{d}
\mid _{Y \times Y} = \widehat{d}_{Y} ; $ and for any $ x \in X $
and $ y \in Y, $ define $ \widehat{d} ( x, y ) = \widehat{d} ( y,
x) = D. $ Then it is easy to see that $ \widehat{d} $ is an
admissible non-Archimedean metric on $ Z. $ So
\begin{align*}
\widehat{d}_{GH} ( X, \ Y ) & = \overline{\widehat{d}}_{GH}( X, \
Y ) \leq \widehat{d}_{H} (X, \ Y ) \\
& = \max \{ \sup _{x \in X } \text{dist} (x, \ Y ), \ \sup _{y \in
Y }\text{dist} (y, \ X ) \} \\
& = \max \{ D, \ D \} = D.
\end{align*}
This proves (3A).
\par  \vskip 0.15 cm
(3B) \ By (3A), we have $ \widehat{d}_{GH} ( X, \ Y ) \leq D. $ We
may assume that $ D = \text{diam}(X) > \text{diam}(Y). $ If $
\widehat{d}_{GH} ( X, \ Y ) < D, $ then we can take an $
\varepsilon $ such that  $$ \max \{ \text{diam}(Y), \
\widehat{d}_{GH} ( X, \ Y ) \} < \varepsilon < D. $$ And then
there exist $ x_{1}, \ x_{2} \in X $ such that $ \widehat{d}_{X} (
x_{1}, \ x_{2} ) > \varepsilon . $ By Theorem 2.23.(1), there
exists a $ strong \ \varepsilon -$isometry $ f : \ X
\longrightarrow Y. $ By the condition $ (SI_{2}) $ of Def.2.22, we
get $ \widehat{d}_{Y} ( f(x_{1}), \ f(x_{2}) ) = \widehat{d}_{X} (
x_{1}, \ x_{2} ) > \varepsilon , $ so $ \text{diam}(Y) \geq
\widehat{d}_{Y} ( f(x_{1}), \ f(x_{2}) ) > \varepsilon . $ A
contradiction! Hence $ \widehat{d}_{GH} ( X, \ Y ) = D. $ This
proves (3B). And the proof of Theorem 4.2 is completed. \quad $
\Box $
\par  \vskip 0.2 cm

{\bf Question 4.3.} \ When will the equality of Theorem 4.2.(1)
hold ?
\par  \vskip 0.2 cm

{\bf Theorem 4.4.} \ Let $ \{ X_{n} \}_{n = 1 }^{\infty } $ be a
sequence of non-Archimedean metric spaces with $ \text{diam} (
X_{n} ) < + \infty $ for each $ n \in \N , $ and $ X $ be a
non-Archimedean metric space with $ \text{diam} ( X ) < + \infty .
$ If $ X_{n} \longrightarrow _{\text{GH}_{S}} X, $ then
\par  \vskip 0.15 cm
(1) \ If $ \text{diam} ( X ) = 0 , $ then $ \text{diam} ( X_{n} )
\longrightarrow 0 $ as $ n \longrightarrow \infty . $
\par  \vskip 0.15 cm
(2) \ If $ \text{diam} ( X ) > 0 , $ then there exists a $ n_{0}
\in \N $ such that for all $ n > n_{0}, \ \text{diam} ( X_{n} ) =
\text{diam} ( X) . $
\par  \vskip 0.2 cm

{\bf Proof.} \ (1) \ If the conclusion that $ \text{diam} ( X_{n}
) \longrightarrow 0 $ as $ n \longrightarrow \infty $ does not
hold, then there exists an $ \varepsilon _{0} > 0 $ such that for
any $ n \in \N , $ there exists a $ N \in \N $ such that $ N > n $
and $ \text{diam} ( X_{N} ) \geq \varepsilon _{0} . $ \\
Since $ X_{n} \longrightarrow _{\text{GH}_{S}} X, $ for $
\varepsilon = \varepsilon _{0} / 2 , $ there exists a $ n_{0} \in
\N $ such that for all $ n > n_{0} , \ \widehat{d}_{GH} ( X_{n}, \
X ) < \varepsilon _{0} / 2 = \varepsilon . $ As discussed above,
there exists $ N \in \N $ with $ N > n_{0} $ such that $
\text{diam} ( X_{N}) \geq \varepsilon _{0} > 0. $ But $
\text{diam} ( X ) = 0 , $ so by Theorem 4.2.(3B), \ $
\widehat{d}_{GH} ( X_{N}, \ X ) = \text{diam} ( X_{N}) \geq
\varepsilon _{0} > \varepsilon. $ a contradiction! Hence $
\text{diam} ( X_{n} ) \longrightarrow 0 $ as $ n \longrightarrow
\infty . $ This proves (1).
\par  \vskip 0.15 cm
(2) \ Take $ \varepsilon = \frac{1}{2}\text{diam} ( X ) > 0, $
then there exists a $ n_{0} \in \N $ such that for all $ n >
n_{0}, \ \widehat{d}_{GH} ( X_{n}, \ X ) < \varepsilon , $ and
then $ \widehat{d}_{GH} ( X_{n}, \ X ) < \max \{ \text{diam} (
X_{n}), \ \text{diam} ( X ) \} $ for all $ n > n_{0}. $ So by
Theorem 4.2.(3B), we get $ \text{diam} ( X_{n}) = \text{diam} ( X
) $ for all $ n > n_{0}. $ This proves (2). And the proof of
Theorem 4.4 is completed. \quad $ \Box $
\par  \vskip 0.2 cm

{\bf Remark 4.5.} \ For a sequence $ \{ X_{n} \}_{n = 1 }^{\infty
} $ of non-Archimedean metric spaces with all $ \text{diam} (
X_{n}) < + \infty . $ By Theorem 4.4, we know that, if $ \{ X_{n}
\}_{n = 1 }^{\infty } $ is a Cauchy sequence under the
non-Archimedean Gromov-Hausdorff metric $ \widehat{d}_{GH}, $ then
either $ \text{diam} ( X_{n}) \longrightarrow 0 $ as $ n
\longrightarrow \infty $ or there exists a $ N \in \N $ such that
$ \text{diam} ( X_{N}) = \text{diam} ( X_{N + i}) $ for all $ i =
1, \ 2, \ \cdots . $
\par  \vskip 0.2 cm

{\bf Corollary 4.6.} \ Let $ X $ be a non-Archimedean metric space
with $ \text{diam} ( X ) < + \infty . $ Let $ \varepsilon > 0 $
and $ Y $ be an $ \varepsilon -$net in $ X . $ If $ \varepsilon <
\text{diam} ( X ) , $ then $ \text{diam} ( Y ) = \text{diam} ( X
). $
\par  \vskip 0.2 cm

{\bf Proof.} \ If $ \text{diam} ( Y ) \neq \text{diam} ( X ), $
then $ \text{diam} ( Y ) < \text{diam} ( X ), $ so by Theorem
4.2.(3B), $ \widehat{d}_{GH} ( X, \ Y ) = \text{diam} ( X ). $ But
by Example 2.2 we know that $ \widehat{d}_{GH} ( X, \ Y ) \leq
\varepsilon . $ A contradiction! \quad $ \Box $
\par  \vskip 0.2 cm

The following are some examples on computing $ \widehat{d}_{GH}
(X, \ Y ) $ in the local fields by using the above tools.
\par  \vskip 0.2 cm
To begin with, as before, let $ p $ and $ q $ be two distinct
rational prime numbers, $ \C_{p} $ and  $ \C_{q} $ be the
corresponding Tate fields endowed with the non-Archimedean metrics
$ \widehat{d}_{p} $ and $ \widehat{d}_{q} $ respectively, where $
\widehat{d}_{p} $ and $ \widehat{d}_{q} $ are induced by the
corresponding normalized non-Archimedean absolute values $ \mid
\cdot \mid _{p} $ and $ \mid \cdot \mid _{q} $ ( i.e., \ $ \mid p
\mid _{p} = \frac{1}{p} $ and \ $ \mid q \mid _{q}
= \frac{1}{q} $ ). \\
We fix an embedding of $ \Q_{p} $ into $ \C_{p}, $ similarly for $
\Q_{q} $ into $ \C_{q}. $ Let $ F \subset \C_{p} $ and $ K \subset
\C_{q} $ be local fields with $ [F : \Q_{p} ] = m $ and $ [K :
\Q_{q} ] = n, $ respectively. We denote \\
$ \mathcal{O}_{F} $ \ the ring of integers of $ F ; \quad U_{F} =
\mathcal{O}_{F}^{\ast } $ \ the unit group of $
\mathcal{O}_{F}; $ \\
$ \mathcal{M}_{F} = \pi _{F} \mathcal{O}_{F} $ \ the maximal ideal
of $ \mathcal{O}_{F} $ with the uniformizer $ \pi _{F} $ for $
\mathcal{O}_{F}; $ \\
$ k_{F} =  \mathcal{O}_{F} / \mathcal{M}_{F} $ \ the residue
field; \quad $ e_{F} $ and $ f_{F} $ be the ramification index and
residue degree, respectively. \ We have $ e_{F} \cdot f_{F} = m .
$ \\
Likewise, for $ K, $ we denote \\
$ \mathcal{O}_{K} $ \ the ring of integers of $ K ; \quad U_{K} =
\mathcal{O}_{K}^{\ast } $ \ the unit group of $
\mathcal{O}_{K}; $ \\
$ \mathcal{M}_{K} = \pi _{F} \mathcal{O}_{K} $ \ the maximal ideal
of $ \mathcal{O}_{K} $ with the uniformizer $ \pi _{K} $ for $
\mathcal{O}_{K}; $ \\
$ k_{K} =  \mathcal{O}_{K} / \mathcal{M}_{K} $ \ the residue
field; \quad $ e_{K} $ and $ f_{K} $ be the ramification index and
residue degree, respectively. \ We have $ e_{K} \cdot f_{K} = n .
$ \\
Obviously, we have $ \mid \pi _{F} \mid _{p} = p^{- 1 / e_{F} } $
\ and \ $ \mid \pi _{K} \mid _{q} = q^{- 1 / e_{K} }. $ \\
( See [K], [L] and [Se] for these and relating facts. )
\par  \vskip 0.2 cm

{\bf Example 4.7.} \ (1) \ $ \widehat{d}_{GH} ( F, \ K ) = \infty
. $ In fact, for any fields $ F^{\prime} \subset \C_{p} $ and $
K^{\prime} \subset \C_{q}, $ we have $ \widehat{d}_{GH} (
F^{\prime}, \ K^{\prime} ) = \infty. $
\par  \vskip 0.15 cm
(2) \ $ \widehat{d}_{GH} ( \mathcal{O}_{F}, \ \mathcal{O}_{K} ) =
1. $
\par  \vskip 0.15 cm
(3) \ For any integers $ s, \ t \in \Z, \quad \widehat{d}_{GH} (
\mathcal{M}_{F}^{s} , \ \mathcal{M}_{K}^{t} ) = \max \{ p^{- s /
e_{F} }, \ q^{- t / e_{K} } \}. $
\par  \vskip 0.2 cm

{\bf Proof.} \ (1) \ If $ \widehat{d}_{GH} ( F, \ K ) < \infty , $
then $ \widehat{d}_{GH} ( F, \ K ) < \varepsilon $ for some $
\varepsilon > 1. $ So by Theorem 2.23.(1), there exists a $ strong
\ \varepsilon -$isometry $ f : \ F \longrightarrow K. $ Since $
\mid F \mid _{p} = p^{\Z / e_{F} } \cup \{ 0 \} , $ there exists $
x_{1}, \ x_{2} \in F $ such that $$ \widehat{d}_{p} ( x_{1}, \
x_{2} ) = \ \mid x_{1} - x_{2} \mid _{p} \ = ( p^{- 1 / e_{F}
})^{\text{ord}_{\pi _{F}}( x_{1} - x_{2})} \geq \varepsilon . $$
Denote $ \text{ord}_{\pi _{F}}( x_{1} - x_{2}) = - r \in \Z, $
then $ r \in \Z _{ > 0} $ and $ \widehat{d}_{p} ( x_{1}, \ x_{2} )
= p^{ r / e_{F} } \geq \varepsilon . $ Then by the condition $
(SI_{2}) $ of Def.2.22, we get $ \widehat{d}_{q} ( f( x_{1}), \ f
( x_{2}) ) = \widehat{d}_{p} ( x_{1}, \ x_{2} ) = p^{ r / e_{F} }.
 $ Write $ y_{i} = f( x_{i}) \in K \ ( i = 1, \ 2 ). $ Then
 $ \text{ord}_{\pi _{K}}( y_{1} - y_{2}) = - r^{\prime } \in \Z,
 $ so $ \widehat{d}_{q} ( f( x_{1}), \ f
( x_{2}) ) = ( q^{- 1 / e_{K} })^{\text{ord}_{\pi _{K}}( y_{1} -
y_{2})} = ( q^{- 1 / e_{K} })^{ - r^{\prime } } = q^{ r^{\prime} /
e_{K} }. $ Thus we get $ 1 < p^{ r / e_{F} } = q^{ r^{\prime} /
e_{K} }. $ Hence $ 1 < p^{ r \cdot e_{K} } = q^{ r^{\prime} \cdot
e_{F} }. $ By the unique factorization of $ \Z, $ this is
impossible because $ p \neq q. $ Therefore $ \widehat{d}_{GH} ( F,
\ K ) = \infty . $ This proves (1).
\par  \vskip 0.15 cm
(2) \ Firstly, since $ \text{diam}(\mathcal{O}_{F}) = \text{diam}
\mathcal{O}_{K} = 1, $ by Theorem 4.2.(3A), we have $
\widehat{d}_{GH} ( \mathcal{O}_{F}, \ \mathcal{O}_{K} ) \leq 1. $
\\
Next we need to prove $ \widehat{d}_{GH} ( \mathcal{O}_{F}, \
\mathcal{O}_{K} ) \geq 1. $ To see this, let $ f : \
\mathcal{O}_{F}  \longrightarrow \mathcal{O}_{K} $ be any $ strong
\ \varepsilon -$isometry with $ \varepsilon > 0. $ We come to
prove that $ \varepsilon > 1. $ Firstly, for the residue fields $
k_{F} $ and $ k_{K} $ as above, we have $ \sharp k_{F} = p^{f_{F}}
$ and $ \sharp k_{K} = q^{f_{K}}. $ Obviously $ \sharp k_{F} \neq
\sharp k_{K} $ because $ p \neq q. $ So we may as well assume that
$ \sharp k_{F} > \sharp k_{K}. $ Take $ A = \{ a_{1}, \ a_{2}, \
\cdots , \ a_{p^{f_{F}}} \} \subset \mathcal{O}_{F} $ such that $
k_{F} = \{ \overline{a_{1}}, \ \overline{a_{2}}, \ \cdots , \
\overline{a_{p^{f_{F}}}} \} $ with $ \overline{a_{i}} = a_{i} \
\text{mod} \mathcal{M}_{F} $ for each $ i. $ Then for the map $ f
$ above, since $ \sharp k_{F} > \sharp k_{K}, $ there exist $
a_{i}, \ a_{j} \in \mathcal{O}_{F} $ with $ i \neq j $ such that $
f( a_{i} ) \equiv f( a_{j} ) \text{mod} \mathcal{M}_{K}. $ If $
\varepsilon \leq 1, $ then note that $ \overline{a_{i}} \neq
\overline{a_{j}}, $ we have $$ \widehat{d}_{p} ( a_{i}, \ a_{j} )
= \mid a_{i} - a_{j} \mid _{p} = ( p^{- 1 / e_{F}
})^{\text{ord}_{\pi _{F}}( a_{i} - a_{j})} = 1 \geq \varepsilon ,
$$ so by the condition $ (SI_{2}) $ of Def.2.22, we have
\begin{align*}
&\widehat{d}_{q} ( f(a_{i}), \ f(a_{j}) ) = \widehat{d}_{p} (
a_{i}, \ a_{j} ) = 1. \quad \text{But} \\
&\widehat{d}_{q} ( f(a_{i}), \ f(a_{j}) ) = ( q^{- 1 / e_{K}
})^{\text{ord}_{\pi _{K}}( f(a_{i}) - f(a_{j}))} \leq q^{- 1 /
e_{K} } < 1.
\end{align*}
A contradiction ! So we must have $ \varepsilon > 1. $ Then by
Corollary 2.24, we get $ \widehat{d}_{GH} ( \mathcal{O}_{F}, \
\mathcal{O}_{K} ) \geq 1. $ Therefore, $ \widehat{d}_{GH} (
\mathcal{O}_{F}, \ \mathcal{O}_{K} ) = 1. $ This proves (2).
\par  \vskip 0.15 cm
(3) \ If $ s = t = 0, $ then this is the above case (2). The other
cases follow directly from Theorem 4.2.(3B) by the fact that $
\text{diam}( \mathcal{M}_{F}^{s} ) = p^{- s / e_{F} } \neq q^{- t
/ e_{K} } = \text{diam} ( \mathcal{M}_{K}^{t} ) ). $ \quad $ \Box
$
\par  \vskip 0.2 cm
{ \bf Note added for Example 4.7.} \ Similarly one can work out
more other examples by the method used above.
\par  \vskip 0.2 cm
Now we come to answer the above Question 2.17.
\par  \vskip 0.2 cm

{\bf Theorem 4.8.} \ The metric ratio function $ \Upsilon _{m} $
is unbounded, in other words, for any $ c \geq 2 , $ there exist
non-Archimedean metric spaces $ X $ and $ Y $ such that $
\widehat{d}_{GH } ( X , Y ) \geq c \cdot d_{GH } ( X , Y ) . $
\par  \vskip 0.2 cm

{\bf Proof.} \ We need to construct a series of such metric spaces
$ X $ and $ Y. $ \\
For $ p-$adic integer ring $ \Z_{p} $ and $ q-$adic integer ring $
\Z_{q} $ with rational primes $ p > q $ as before. We first
construct the following set $$ \Z_{q}^{\Delta } = \Z_{q} \sqcup \{
t_{q}, \ \cdots , \ t_{p - 1} \} \ ( \text{ the disjoint union })
$$ with the $ p - q $ number of indeterminate elements
$ t_{q}, \ \cdots , \ t_{p - 1}. $ Then we define a function $
\widehat{d}_{\Delta } $ on $ \Z_{q}^{\Delta } \times
\Z_{q}^{\Delta } $ as follows: \\
$ \widehat{d}_{\Delta } \mid _{Z_{q} \times Z_{q} } =
\widehat{d}_{q}; \quad \widehat{d}_{\Delta } ( t_{i}, \ t_{i}) = 0
\quad ( \forall \ i \in \{ q, \ \cdots , \ p - 1 \});  \\
\widehat{d}_{\Delta } ( t_{i}, \ t_{j}) = 1 + \frac{1}{q} \quad (
\forall \
i, \ j \in \{ q, \ \cdots , \ p - 1 \} $ and $ i \neq j ); $ \\
$ \widehat{d}_{\Delta } ( t_{i}, \ a ) = \widehat{d}_{ \Delta } (
a, \ t_{i} ) = 1 + \frac{1}{q} \quad (
\forall \ i \in \{ q, \ \cdots , \ p - 1 \} $ and $ a \in \Z_{q}). $ \\
It is easy to verify that $ \widehat{d}_{\Delta } $ is a
non-Archimedean
metric on $ \Z_{q}^{\Delta }. $ \\
Since $ \text{diam} ( \Z_{p}) = 1 $ and $ \text{diam} (
\Z_{q}^{\Delta } ) = 1 + \frac{1}{q}, $ by Theorem 4.2.(3B), we
get $$ \widehat{d}_{GH } ( \Z_{p}, \ \Z_{q}^{\Delta } ) = 1 +
\frac{1}{q}. $$ Next we come to prove that $$ d_{GH } ( \Z_{p}, \
\Z_{q}^{\Delta } ) \leq  \frac{1}{2 q }. $$ To see this, we define
a correspondence $ \mathcal{C} $ between $ \Z_{p} $ and $
\Z_{q}^{\Delta } $ as follows: \\
Firstly, we have
\begin{align*}
&\Z_{p} = \sqcup _{i = 0 }^{p - 1 } A_{i} \ ( \text{ the disjoint
union }) \ \text{with} \ A_{i} = i + p \Z_{p} \\
&\Z_{q}^{\Delta } = ( \sqcup _{i = 0 }^{q - 1 } B_{i} ) \sqcup (
\sqcup _{i = q }^{p - 1 } \{ t_{i} \} ) \ \text{with} \ B_{i} = i
+ q \Z_{q}.
\end{align*}
Then we set $$ \mathcal{C} =  \cup _{i = 0 }^{q - 1 } ( A_{i}
\times  B_{i} ) \cup \cup _{i = q }^{p - 1 } ( A_{i} \times \{
t_{i} \} ) \subset \Z_{p} \times \Z_{q}^{\Delta }. $$ Obviously $
\mathcal{C} $ is a correspondence between $ \Z_{p} $ and $
\Z_{q}^{\Delta } . $
\par  \vskip 0.15 cm
Assertion. \ we have $ \text{dis} \mathcal{C} \leq 1 / q . $
\par  \vskip 0.15 cm
To see this, for any $ ( x_{1}, \ y_{1} ), \ ( x_{2}, \ y_{2} )
\in \mathcal{C}, $ we discuss into the following cases: \\
(a) \ $ x_{1}, \ x_{2} \in A_{i} . $  \\
(a$_{1}$) \ If $ 0 \leq i \leq q - 1,  $ then $ y_{1}, \ y_{2} \in
B_{i}, $ so $$ \mid \widehat{d}_{p} ( x_{1}, \ x_{2} ) -
\widehat{d}_{q} ( y_{1}, \ y_{2} ) \mid \ \leq \max \{ \text{diam}
( A_{i}), \ \text{diam} ( B_{i}) \} = \max \{ \frac{1}{p}, \
\frac{1}{q} \} = \frac{1}{q}. $$ \\
(a$_{2}$) \ If $ q \leq i \leq p - 1,  $ then $ y_{1} = y_{2} =
t_{i}, $ so \ $ \mid \widehat{d}_{p} ( x_{1}, \ x_{2} ) -
\widehat{d}_{\Delta } ( y_{1}, \ y_{2} ) \mid \ = \ \mid
\widehat{d}_{p} ( x_{1}, \ x_{2} ) \mid \ \leq \text{diam} (
A_{i}) = 1 / p < 1 /
q. $ \\
(b) \ $ x_{1} \in A_{i} $ and $  x_{2} \in A_{j} \ ( i \neq j ). $
\\
(b$_{1}$) \ If $ i, \ j \leq q - 1, $ then $ y_{1} \in B_{i} $ and
$ y_{2} \in B_{j}. $ So $ \mid \widehat{d}_{p} ( x_{1}, \ x_{2} )
- \widehat{d}_{\Delta } ( y_{1}, \ y_{2} ) \mid \ = \ \mid  1 - 1
\mid \ = 0. $ \\
(b$_{2}$) \ If $ i, \ j > q - 1, $ then $ y_{1} = t_{i} $ and $
y_{2} = t_{j}. $ So \\
$ \mid \widehat{d}_{p} ( x_{1}, \ x_{2} ) - \widehat{d}_{\Delta }
( y_{1}, \ y_{2} ) \mid \ = \ \mid  1 - \widehat{d}_{\Delta } (
t_{i}, \ t_{j} ) \mid \ = \ \mid  1 - ( 1 + \frac{1}{q}) \mid \ =
1 / q. $ \\
(b$_{3}$) \ (The other cases ) \ We may as well assume that $ i
\leq q - 1 $ and $ j > q - 1, $ then $ y_{1} \in B_{i} $ and $
y_{2} = t_{j}. $ So \\
$ \mid \widehat{d}_{p} ( x_{1}, \ x_{2} ) - \widehat{d}_{\Delta }
( y_{1}, \ y_{2} ) \mid \ = \ \mid  1 - ( 1 + \frac{1}{q}) \mid \
= 1 / q. $ \\
To sum up, we obtain that $$ \text{dis} \mathcal{C} = \sup \{ \mid
\widehat{d}_{p} ( x_{1}, \ x_{2} ) - \widehat{d}_{\Delta } (
y_{1}, \ y_{2} ) \mid \ : \ ( x_{1}, \ y_{1} ), \ ( x_{2}, \ y_{2}
) \in \mathcal{C} \} \leq \frac{1}{q}. $$ This proves the above
assertion. So by Theorem 7.3.25 of [BBI, p.257 ], we get $$ d_{GH}
( \Z_{p}, \ \Z_{q}^{\Delta } ) \leq \ \frac{1}{2} \ \text{dis}
\mathcal{C} \leq \frac{1}{2 q }. \quad \text{Therefore,} $$ \ $$
\Upsilon _{m} ( \Z_{p}, \ \Z_{q}^{\Delta } ) = \frac{
\widehat{d}_{GH} ( \Z_{p}, \ \Z_{q}^{\Delta } ) }{ d_{GH} (
\Z_{p}, \ \Z_{q}^{\Delta } ) } \geq \frac{ 1 + \frac{1}{q}}{
\frac{1}{2 q }} = 2 q + 2 \longrightarrow \infty \ \text{as} \ q
\longrightarrow \infty . $$ This proves Theorem 4.8. \quad $ \Box
$
\par  \vskip 0.2 cm

{\bf Remark 4.9.} \ It is well known that $ d_{GH} ( X, \ Y ) \geq
\frac{1}{2} \ \mid \text{diam} (X) - \text{diam} (Y) \mid $ \ for
any metric spaces $ X $ and $ Y $ with $ \text{diam} (X) < \infty
$ \ ( see Exercise 7.3.14 of [BBI, p.255] ), so we have $$ d_{GH}
( \Z_{p}, \ \Z_{q}^{\Delta } ) \geq \frac{1}{2} \ \mid \text{diam}
( \Z_{p} ) - \text{diam} ( \Z_{q}^{\Delta } ) \mid \ = \frac{1}{2}
\ \mid \ 1 - ( 1 + \frac{1}{q} ) \mid \ = \frac{1}{2 q }. $$ Hence
by the computation result in the above Theorem 4.8, we obtain that
$$ d_{GH}( \Z_{p}, \ \Z_{q}^{\Delta } ) = \frac{1}{2 q }. \quad
\text{Therefore,}  $$  \
$$ \Upsilon _{m} ( \Z_{p}, \ \Z_{q}^{\Delta } )
= \frac{ \widehat{d}_{GH} ( \Z_{p}, \ \Z_{q}^{\Delta } ) }{ d_{GH}
( \Z_{p}, \ \Z_{q}^{\Delta } ) } = \frac{ 1 + \frac{1}{q}}{
\frac{1}{2 q }} = 2 q + 2.  \quad  \Box  $$
\par  \vskip 1 cm
\hspace{-0.8cm} {\bf References }
\begin{description}

\item[[BBI]] D. Burago, Y. Burago, S. Ivanov, A Course in Metric
Geometry, Providence, Rhode Island: American Mathematical Society,
2001.

\item[[BGR]] S. Bosch, U. Guntzer, R. Remmert, Non-Archimedean
Analysis, Berlin: Springer-Verlag, 1984.

\item[[G]] M. Gromov, Metric Structures for Riemannian and
Non-Riemannian Spaces, Boston: Birkhauser, 2001.

\item[[K]] N. Koblitz, $ p-$adic Numbers, $ p-$adic Analysis, and
Zeta Functions, Second Edition, New York: Springer-Verlag, 1984.

\item[[L]] S. Lang, Algebraic Number Theory, Second Edition, New
York: Springer-Verlag, 1994.

\item[[Q1]] D. R. Qiu, Arithmetic applications of Non-Archimedean
Gromov-Hausdorff metric, preparation.

\item[[Q2]] D. R. Qiu, Metric structures on $ p-$adic manifolds,
preparation.

\item[[Sc]] W. H. Schikhof, Ultrametric Calculus, London:
Cambridge University Press, 1984.

\item[[Se]] J.-P. Serre, Local Fields, New York: Springer-Verlag,
1979.

\item[[Z]] I. Zarichnyi, Gromov-Hausdorff ultrametric { \it arXiv:
math} / 0511437v1, 17 Nov. 2005.

\end{description}

\end{document}